\input amstex.tex

\input amsppt.sty

\TagsAsMath

\NoRunningHeads \magnification=1200

\hsize=5.0in\vsize=7.0in

\hoffset=0.2in\voffset=0cm

\nonstopmode

\document
\topmatter

\title{On instability of excited states of the nonlinear Schr\"odinger equation
}
\endtitle

\author
Scipio Cuccagna
\endauthor

\address
DISMI University of Modena and Reggio Emilia, via Amendola 2,
Padiglione Morselli, Reggio Emilia 42100 Italy
\endaddress

\email cuccagna.scipio\@unimore.it\endemail

\abstract We introduce a new notion of linear stability for standing
waves of the nonlinear Schr\"odinger equation  (NLS) which  requires
not only that the spectrum of the linearization  be real, but also
that the generalized kernel    be not degenerate and that the
signature of all the positive eigenvalues be positive. We prove that
excited states of the NLS are not linearly stable in this more
restrictive sense. We then give a partial proof that this more
restrictive notion of linear stability is a necessary condition to
have orbital stability.

\endabstract

\endtopmatter

\head \S 1 Introduction \endhead

We  consider the NLS
$$iu_t +\Delta u +\beta ( |u|^2 )u=0\, ,\, (t,x)
\in \Bbb R \times \Bbb R^3 \, , \,  u(0,x)=u_0(x) . \tag 1.1$$ We
will assume the following hypotheses.

{\item {(H1)}} There exists an open interval $\Cal{O}\subseteq
(0,+\infty )$   such that $
    \Delta u-\omega u+\beta(u^2)u=0
$ admits  a family of standing waves $e^{it \omega }\phi _ {\omega }
(x)$ with $\phi _ {\omega }$ real valued  and   $\omega \in \Cal
O\to \phi _ {\omega } \in  C^1(\Cal{O}, H^1_{r }(\Bbb R^3) )$, with
$H^1_{r }$ denoting radially symmetric finite energy functions.

{\item {(H2)}} $\beta (0)=0$, $\beta\in C^1(\Bbb R,\Bbb R)$.

{\item {(H3)}} There exists a $1<p< 5 $ such that for every $k=0,1$,
$$\left| \frac{d^k}{dv^k}\beta(v^2)\right|\lesssim
|v|^{p-k-1} \quad\text{if $|v|\ge 1$}.$$

\medskip
In this paper we focus on the question of orbital stability of the
excited standing waves $e^{it \omega }\phi _ {\omega }  $. This
question has been explored thoroughly  for ground states, i.e. when
we can pick $\phi _ {\omega } (x)> 0$ for all $x$, see
\cite{CL,We1,GSS1-2} and countless other papers. In the case of
ground states, orbital stability generally is equivalent to the
standard notion of linear stability, which we review now. Recall
that for $\sigma _j$ the Pauli matrices  given below, then the
linearization $H_\omega$ is given (see \S 2)  by
$$\aligned &\sigma _1=\left[ \matrix  0 &
1  \\
1 & 0
 \endmatrix \right] \, ,
\sigma _2=\left[ \matrix  0 &
i  \\
-i & 0
 \endmatrix \right] \, ,
\sigma _3=\left[ \matrix  1 &
0  \\
0 & -1
 \endmatrix \right]
\, , \\&H_{\omega }=\sigma _3 \left [ - \Delta + \omega  -  \beta
(\phi ^2 _{\omega  }) - \beta ^\prime (\phi ^2 _{\omega  })\phi ^2
_{\omega   } \right ] +i  \beta ^\prime (\phi ^2 _{\omega  })\phi ^2
_{\omega  }
 \sigma _2    .
\endaligned \tag 1.2$$
A ground state is generally orbitally stable if $\sigma (H_\omega )
\subset \Bbb R$, but not always    since   there are well known
cases like the critical equation $iu_t+\Delta
u+|u|^{\frac{4}{n}}u=0$ with $\sigma (H_\omega ) \subset \Bbb R$
where the generalized kernel of $H_\omega $ has higher degeneracy
and the ground state is orbitally unstable.   See also the work by
Comech and Pelinovsky \cite{CoP}.  In this paper we look at excited
states. More precisely   assume that the  $\phi _ {\omega } (x) $
are real valued  and change sign. For these standing waves less
appears to be known. One can  look at the spectrum $\sigma (H_\omega
)$ of the linearization (1.2), which is symmetric with respect of
the coordinate axes. It has been known for a long time, but proved
rigorously only recently by Mizumachi \cite{M1} in dimension $ 2$
(the argument
  extends to any dimension), that if $\sigma (H_\omega )\not
\subset \Bbb R$ then $\phi _ {\omega }  $ is not orbitally stable in
$H^1 (\Bbb R^n)$. In the literature various studies of instability
of standing waves are based on this form of linear instability, see
\cite{J,Gr1-2,M2-4}.  Based on the above considerations, classically
a standing wave is called linearly stable if $\sigma (H_\omega )
\subset \Bbb R$. While this classical notion of linear stability is
adequate in the case of ground states, our purpose here is to show
that  it is inadequate in the case of excited states, and to propose
a substitute. In \S 2 Definition 2.3 we give a new definition  of
linear stability. Succinctly, this requires not only $\sigma
(H_\omega )\subset \Bbb R$, but also that the generalized kernel
$N_g(H_\omega )$  be not degenerate and that the signature of all
positive eigenvalues    be positive.  This definition of linear
stability generally coincides with the classical one in the case of
ground states, because in that case $H_\omega$ has no positive
eigenvalues of negative signature. However for excited states we
have:

\proclaim{Theorem 1.1} Consider hypotheses (H1-3) and suppose that
the  $\phi _ {\omega }  $ are real valued and change sign. Then
$\phi _ {\omega }  $ is not linearly stable in the sense of
Definition 2.3.
\endproclaim

 The fact that
excited states do not meet a new and more stringent definition of
linear stability is by itself  not very significant. What matters is
to see whether this new definition sheds some light on the question
of orbital instability of excited states. In this respect   we
conjecture that a standing wave $e^{it\omega }\phi _ {\omega } $
satisfying (H1-3) and with  $\phi _ {\omega } $ real valued is
orbitally stable if and only if it satisfies Definition 2.3 (we also
conjecture that in that case it is also asymptotically stable in the
sense of \cite{CM}). In sections 3 and 4 we establish in special
situations part of the conjecture, that is that, under appropriate
hypotheses, excited states with $\sigma (H_\omega )\subset \Bbb R$
are nonetheless orbitally unstable.

For an excited state with $\sigma (H_\omega )\subset \Bbb R$ there
are three mechanisms which   yield orbital instability, two known
and a third one explored here. The first two mechanisms are
basically linear, because in these two cases, even though $\sigma
(H_\omega )\subset \Bbb R$, there are arbitrarily small
perturbations of $H_\omega $, of  appropriate restricted type, with
eigenvalues outside $\Bbb R$. These   first two mechanisms are also
non generic. The third mechanism, the only one explored here, is
nonlinear   and applies to cases where the condition $\sigma
(H_\omega )\subset \Bbb R$ is stable under perturbation. The first
mechanism of instability arises from the degeneracy of the
generalized kernel $N_g(H_\omega )$. This first mechanism is
explored  in \cite{CoP} and is not discussed here. So in sections 3
and 4 we assume that $N_g(H_\omega )$ is non degenerate, which is a
generic condition.   The second mechanism of instability  is related
to the possible presence of eigenvalues of negative signature
embedded in $\sigma _c(H_\omega ) $. This   phenomenon is absent for
ground states. While we cannot point to examples  in the literature
of this occurrence for excited states,
 it should be possible to prove their existence   via
perturbation theory such as \cite{CHM}. Then   orbital instability
should follow  by   essentially linear mechanisms of the type in
\cite{G1,TY4,CPV}. We do not discuss the above two mechanisms and we
only say that if   $N_g(H_\omega )$ is degenerate and  there are
embedded eigenvalues of negative signature, there are additional
reasons for instability with respect to the ones described here.
Furthermore, if present, the first two mechanisms   will usually
prevail, because usually they unfold more rapidly than the third.
The third mechanism is nonlinear and robust. The setting is related
to attempts in a long list of papers
 \cite{BP,SW2,TY1-3,Cu3,BS,T,GS,SW1,CM,Cu1,CT} to prove asymptotic
 stability of stable ground states. We assume
 more regularity   on the nonlinearity $\beta (r)$.  This because we consider appropriate
 Taylor expansions of  $\beta (|w|^2)w$ and   normal forms
 transformations which lead  to a particular expansion of equation (1.1)
 around the orbit of an excited state. There is a natural
 decomposition in discrete and continuous modes, with the discrete
 ones satisfying a perturbation of a Hamiltonian system. In sections
3 and 4
  it is described, after \cite{BP,SW2,TY1-3,Cu3,BS,T,SW1,GS,CM}, a possible mechanism
 through which
 the coupling of discrete with the continuous modes  breaks the conservation
 laws of the unperturbed system of discrete modes  and yields, in sections 3 and 4, orbital instability of
 excited states. This mechanism is called Nonlinear Fermi
Golden Rule (FGR), after Sigal \cite{Si}. In section 3 we consider
the case when the portion $\sigma _d(H_\omega )\backslash \{ 0 \}$
of the discrete spectrum is close to the continuous spectrum $\sigma
_c(H_\omega ) $. In this case our proof is valid generically.
 In
section 4,  $\sigma _d(H_\omega )\backslash \{ 0 \}$ is not any more
close to $\sigma _c(H_\omega ) $  and   our proof hinges on a
conjecture on the  FGR, which we assume as hypothesis in Hypothesis
4.4 and is related to similar conjectures in \cite{GS,Gz,CM}. Notice
that even though the conjecture on the  FGR in our present setting
gives orbital
 instability, in other settings, see the papers just referenced, this same mechanism yields
 asymptotic stability.
   The  FGR can be viewed as a
consequence of identities between some coefficients in the system on
continuous and discrete modes. These are Taylor coefficients of the
right hand sides of the equations. If the system derived from a real
valued Hamiltonian $\Cal H$, these coefficients would be mixed
derivatives of $\Cal H$, with different order of differentiation,
and would be equal by the Schwarz lemma on mixed derivatives. Notice
also that the NLS (1.1) is derived from a   real valued Hamiltonian.
Unfortunately we are not able to retain this Hamiltonian structure
during the normal forms argument. Therefore the FGR remains a
conjecture. Another ingredient required is that the above mentioned
coefficients  do not vanish  on appropriate spheres of the
 phase space associated to $H_\omega$, see for instance (3.6). In this respect, we
refer to a question in p.69 \cite{SW2} on   the relation between
vanishing and integrability. At least in the non integrable setting
we hope to have identified the mechanisms of instability. The same
proof up to cosmetic changes can be used for non translation
invariant equations of the form

$$iu_t +( \Delta -q(x))u  +a(x)\beta ( |u|^2 )u=0\tag 1.3$$
for  $q(x)$ short range and   regular and $a(x)$ regular and
bounded. When  $ -\Delta +q(x)$ has two or more negative
eigenvalues, it is easy to manufacture by bifurcation, in the spirit
of \cite{SW1-2,TY1-4,T} for systems of the form (1.3),
 small excited states   for which
 our proof of orbital instability holds. But our proof
 is valid more generally.

In the sequel a  matrix will be called real when its components are
real valued. A matrix will be exponentially decreasing when its
components are functions converging exponentially to 0 as $|x|\to
\infty $. For $f(x)$ and $g(x)$ column vectors, their inner product
is $ \langle f ,
  g\rangle =\int _{\Bbb R^d}{^tf(x)}\cdot  { {\overline{g}(x)}} dx $  with $\overline{g}$ the vector
  with entries which are complex
  conjugate and with $^tf$ the transposed vector. The
  adjoint $H^\ast$ is defined by $
\langle Hf ,
  g\rangle =
\langle f ,
  H^\ast g\rangle .$ By  $L^p(\Bbb R^3)$ or $L^p$ we will
  denote not only the usual $L^p(\Bbb R^3, \Bbb C)$ but also $L^p(\Bbb R^3, \Bbb
  C^2)$, with the exact meaning  clarified by the context. Same
  convention   for the Sobolev spaces $W^{k,p}$, with
  $H^k=W^{k,2}$. For $\langle x \rangle = \sqrt{1+|x|^2} $ we will consider   weighted norms $ \|u\|_{L^{2,s}}=\|\langle  x\rangle ^su\|_{L^2( \Bbb R^3)}
 $,
mixed norms
  $  \| g(t,x) \| _{L^p_tL^q_x}= \|  \| g(t,x) \| _{ L^q_x} \| _{L^p_t
  } $ and $  \| g(t,x) \| _{L^p_tL^{2,s}_x}= \|  \| g(t,x) \| _{ L^{2,s}_x} \| _{L^p_t
  } $. Given a norm $   \| g  \| _{L^p_tL^q_x}$ the pair $(p,q)$
  is said to be admissible if
$\frac{1}{p}=\frac{3}{2}\left ( \frac{1}{2}-\frac{1}{q} \right )$
with $2\le q\le 6$. Given an operator $H$, we set
$R_{H}(z)=(H-z)^{-1}$. In the sequel, for $\lambda \in \Bbb R$ we
will write $ R_{H}^{\pm }(\lambda )=R_{H}(\lambda \pm i0)$ with on
the left an appropriate (i.e. radial or nontangential) limit of
$R_{H}(z)$ for $z\to \lambda $ with $\pm \Im z>0$. Here given a
complex number $z=x+iy$,  we set $\Re z=x$ and $\Im z=y$. For a
matrix or vector $A$, we denote by $^tA$ the transpose.

\head \S 2 Definition of linear stability and proof of Theorem 1.1
\endhead

We can write the following ansatz:
$$\align & u(t,x) = e^{i \vartheta (t)} (\phi _{\omega (t)} (x)+ r(t,x)) \text{ with $\vartheta (t)
=\int _0^t\omega (s) ds + \gamma (t)$ } .\tag 2.1
\endalign$$
Inserting (2.1) in (1.1) we get for some $n(r,\overline{r}
)=O(r^2)$,   $\overline{n(r,\overline{r} )}=n(\overline{r}, {r}  )$
$$\aligned &
  i r_t  =
 -\Delta r +\omega (t) r-
\beta ( \phi _{\omega (t)} ^2 )r -\beta ^\prime ( \phi _{\omega (t)}
^2 )\phi _{\omega (t)} ^2 r \\&-
 \beta ^\prime ( \phi _{\omega (t)} ^2 )
\phi _{\omega (t)} ^2  \overline{  r }+ \dot \gamma (t) \phi
_{\omega (t)} -i\dot \omega (t)
\partial _\omega \phi   _{\omega (t)}
+ \dot \gamma (t) r
 +
  n(r,\overline{r} ) .\endaligned
$$
{Set }  $^tR= ( r
  ,
 \overline{ r}) $, $ ^t\Phi _\omega
=( \phi _{\omega } , \phi _{\omega }) $, $ ^tN(R) =(
n(r,\overline{r} ), - \overline{n} (r,\overline{r} ) ).
 $
We rewrite   the   equation for $r$ as
$$\aligned & i  R _t =H _{\omega (t)}   R +\sigma _3 \dot \gamma (t)  R
+\sigma _3 \dot \gamma (t)\Phi _\omega- i \dot \omega (t)\partial
_\omega \Phi _\omega+N(R ).\\&  H_\omega =\sigma _3 \left [-\Delta +
\omega   -  \beta (\phi ^2 _{\omega  }) - \beta ^\prime (\phi ^2
_{\omega  })\phi ^2 _{\omega  } \right ] +i \sigma _2 \beta ^\prime
(\phi ^2 _{\omega  })\phi ^2 _{\omega  }.\endaligned\tag 2.2 $$ For
an operator $L$ the generalized kernel is the space $N_g(L)= \cup
_{j\ge 1} \ker (L^j)$. We have:

{\item {(a)}}$\sigma _1H_\omega = -H_\omega \sigma _1$ {\item {(b)}}
if we set $H_\omega=\sigma _3(-\Delta +\omega ) +V_\omega (x)$, then
the matrix $V_\omega (x)$ has real valued entries.

\noindent (a)-(b) imply that point spectrum
 $ \sigma _p (   H_\omega )$ and essential spectrum
 $ \sigma _e (  H_\omega )$ are symmetric with respect to  the
 coordinate axes.
The following standard lemma, which requires the exponential decay
of $\phi (x)$ at infinity, is proved in \cite{CPV}:

\proclaim{Lemma 2.1} The point spectrum  is a finite set and we have
 $ \sigma _p (   H_\omega )= \sigma _p (  H_\omega ^\ast )$. Similarly for the
 essential spectrum
 $ \sigma _e (  H_\omega )= \sigma _e (  H_\omega ^\ast
 )=(-\infty , -\omega ]\cup [\omega ,+\infty).$ For each $z\in \sigma _p (  H_\omega )$ the corresponding
generalized eigenspace $N_g(  H_\omega -z) $ has finite dimension.
\endproclaim

We are assuming that $\sigma   (   H_\omega )\subset \Bbb R$ because
otherwise by \cite{M1} the standing wave is unstable. We define:

\proclaim{Definition 2.2} Let $\lambda >0$ be an eigenvalue of
$H_\omega$. We say that $\lambda$ has positive (resp. negative)
signature if the following two points hold:

{\item {(1)}} the algebraic and geometric multiplicity coincide,
i.e. $N_g(H_\omega -\lambda )=\ker (H_\omega -\lambda )$;

 {\item {(2)}} for any $\xi \in \ker (H_\omega -\lambda )$ with $\xi \neq 0$ we
have $\langle \xi , \sigma _3 \xi \rangle >0 $ (resp.  $\langle \xi
, \sigma _3 \xi \rangle <0 $).
\endproclaim
{\it Remark}. Notice that if $z\in \sigma _p  (   H_\omega )
\backslash\Bbb R $ then $\langle \xi , \sigma _3 \xi \rangle =0$ for
any $\xi \in \ker (H_\omega -z )$.
\medskip
{\it Remark}. We are unable to reference  examples of  eigenvalues
$\lambda
>\omega $ with negative signature so we sketch what seems a natural
way to manufacture them. Consider a short range Schr\"odinger
operator $h=-\Delta +q(x)$ with $\sigma _d(h)\supseteq \{
-E_0,-E_1\}$, with $-E_0$ the smallest eigenvalue and with
$-E_0<-E_1<0$. By bifurcation, equation (1.3) with $a(x)\equiv 1$
will have small amplitude excited states $e^{i\omega t}\phi _\omega
(x)$ where $\omega \in \Omega$, $\Omega$ a small open interval with
$E_1$ an endpoint. $\sigma _3(h+\omega )$ will have a real  negative
signature eigenvalue $E_0-\omega$. This will be in $[\omega ,\infty
)$ if $E_0-E_1>E_1$. We have that $H_\omega = \sigma _3(h+\omega
)+V(\phi _\omega )$ with $V(\phi _\omega )(x)=\sigma _3 \varphi
_{\omega}(x)+i\sigma _2\psi _{\omega}(x) $ with $\phi _{\omega}$ and
$\psi _{\omega}$ real valued exponentially decreasing functions.
Generically, by \cite{G1,TY4,CPV}, $H_\omega$ will have a pair of
non real eigenvalues close to $E_0-\omega$. However, in analogy to
the conjecture  in the context of the $N$ body problem in
\cite{AHS}, there is a hypersurface $\Sigma $ of pairs $(\phi
_{\omega},\psi _{\omega})$ such that $H_\omega = \sigma _3(h+\omega
)+\sigma _3 \varphi _{\omega}(x)+i\sigma _2\psi _{\omega}(x) $ has a
real eigenvalue $\lambda (\omega )$ near $E_0-\omega$,   of negative
signature. This is easy to see by employing the Weinstein-Aroszajn
formula in \S 5 \cite{CPV} which relates $\lambda (\omega )$ to
$(\phi _{\omega},\psi _{\omega})$. Furthermore it should be possible
to find some $\beta (|u|^2)u$ in (1.3) such that we have $V(\phi
_\omega )(x)=\sigma _3 \varphi _{\omega}(x)+i\sigma _2\psi
_{\omega}(x) $ with  $(\phi _{\omega},\psi _{\omega})\in \Sigma $.
Notice that the conjecture in \cite{AHS},  proved  for Wigner-Von
Neuman potentials   in \cite{CHM}, is in  a setting   much harder
then ours since we are considering only very    short range
potentials.

\medskip
{\it Remark}. It is not known if there are eigenvalues $\lambda
>\omega $ with positive signature. We conjecture that they do not
exist.

\medskip
{\it Remark}. In \cite{CPV} it is proved that generically
eigenvalues    $\lambda
>\omega $  do not exist in our setting, that is with matrix potentials $V_\omega (x)=\sigma _3 \varphi
_{\omega}(x)+i\sigma _2\psi _{\omega}(x) $ with $\phi _{\omega}$ and
$\psi _{\omega}$ real valued exponentially decreasing functions.

\bigskip

We introduce now our definition of linear stability. The usual
definition is  that $\sigma   (   H_\omega )\subset \Bbb R$. We
prefer the following more stringent definition.

\proclaim{Definition 2.3(Linear stability)}  We will say that $\phi
_\omega $ is linearly stable if the operator $ H_\omega $ satisfies
the following three conditions:
\medskip
{\item {(1)}}  $\sigma   (   H_\omega )\subset \Bbb R$;

   \medskip {\item {(2)}}  if $\lambda >0$ is an eigenvalue of
$H_\omega$ then $\lambda$ has positive signature;

\medskip {\item {(3)}}
     $N_g(H  _\omega )$ is spanned by $ \{ \sigma _3
\Phi _\omega   ,  \partial _\omega \Phi  _\omega     , \partial _x
\Phi _{\omega  }   , \sigma _3x\Phi _\omega     \}   .$

\endproclaim
{\it Remark}. That $N_g(H  _\omega )\supseteq \text{span}\{ \sigma
_3 \Phi _\omega   ,  \partial _\omega \Phi  _\omega     , \partial
_x \Phi _{\omega  }   , \sigma _3x\Phi _\omega     \}$ follows by
direct computation. The fact that  generically  $N_g(H  _\omega )$
strictly larger than this span implies orbital instability has been
explored in \cite{CoP}. In sections 3 and 4 we will assume that, in
the context of even functions, $N_g(H _\omega )=\text{span}\{ \sigma
_3\Phi _\omega ,
\partial _\omega \Phi _\omega \} $.

\medskip {\it Remark}. $N_g(H  _\omega
)=\text{span}\{ \sigma _3 \Phi _\omega   ,  \partial _\omega \Phi
_\omega     , \partial _x \Phi _{\omega  }   , \sigma _3x\Phi
_\omega     \} $  is proved in  \cite{We2}         under hypothesis
(H4), see \S 3, and if $\frac{d}{d\omega }\| \phi _\omega \|
_{2}\neq 0.$

\medskip {\it Remark}. If we break the translation invariance of the
equation, then (3) is replaced by $N_g(H  _\omega )=\text{span}\{
\sigma _3 \Phi _\omega   ,  \partial _\omega \Phi _\omega       \}
$.

\bigskip

{\it Proof of Theorem 1.1.} Set $  L_+ =-\Delta + \omega -\beta
(\phi ^2_\omega   ) -2\phi ^2 _\omega \beta ^\prime (\phi ^2_\omega
) $ and $ L_- =-\Delta + \omega
 -\beta (\phi _\omega  ^2). $  For
$ U= \left[\matrix 1 &
1  \\
i & -i \endmatrix \right]    $ we have $U  \sigma _3H_\omega U^{-1}
= \text{diag} (L_+,L_-)$. Notice that $L_-\phi _\omega =0$. Since
$\phi _\omega$ has  nodes, $L_-$ has a smallest strictly negative
eigenvalue. The corresponding
   ground states of $L_-$ are spherically symmetric.   $L_+\partial
_{x_j}\phi _\omega =0$  for all $j$ and so also $L_+$ has a smallest
strictly negative eigenvalue with corresponding
   ground states which are spherically symmetric. From now on
   in this proof we consider $H_\omega$, $L_+$ and  $L_-$ as acting
   on spherically symmetric functions only. Let
   $$N(\sigma _3H_\omega )=\sum _{\lambda \in \sigma _p(\sigma _3H_\omega )\cap (-\infty , 0)}
\dim \ker (\sigma _3H_\omega -\lambda).$$ Since  we are restricting
to spherically symmetric functions, condition (3) becomes $N_g(H
_\omega )=\text{span}\{ \sigma _3 \Phi _\omega   ,  \partial _\omega
\Phi _\omega       \} $. The fact, discussed
   above, that for both signs  we have
   $\sigma (L_\pm ) \cap (-\infty , 0)\neq \emptyset$,
 implies
$N(\sigma _3H_\omega )\ge 2$. We  claim now that if $H_\omega $
satisfies the conditions in Definition 2.3  we have $$\langle \sigma
_3H_\omega u, u \rangle \ge 0  \text{ for any non zero $u\in H^1\cap
N^\perp _g (H_\omega ^\ast )$}.\tag 2.3$$ Before proving (2.3), we
show that  (2.3) implies $N(\sigma _3H_\omega )\le 1$. We have an
$H_\omega $ invariant splitting $ L^2= N  _g (H_\omega )\oplus
N^\perp _g (H_\omega ^\ast ) $  with $\sigma _3N^\perp _g (H_\omega
^\ast )=N^\perp _g (H_\omega
  )$. This implies that given a generic $u=v+ w\in N  _g (H_\omega
)\oplus N^\perp _g (H_\omega ^\ast ) $,

$$\aligned & \langle
\sigma _3H_\omega u, u \rangle = \langle \sigma _3H_\omega v, v
\rangle  + \langle \sigma _3H_\omega w, w \rangle   \ge \langle
\sigma _3H_\omega v, v \rangle .
\endaligned $$
By (3) Definition 2.3 we have $v=\lambda \sigma _3\Phi +\mu \partial
_\omega \Phi $ with
$$ \aligned & \langle
\sigma _3H_\omega v, v \rangle = \lambda ^2 \langle \sigma
_3H_\omega \sigma _3\Phi,  \sigma _3\Phi \rangle  +  \lambda \mu
\left [ \langle \sigma _3H_\omega \sigma _3\Phi , \partial _\omega
\Phi \rangle + \langle \sigma _3H_\omega \partial _\omega \Phi,
\sigma _3\Phi \rangle \right ] \\& +\mu ^2 \langle \sigma _3H_\omega
\partial _\omega\Phi,  \partial _\omega   \Phi \rangle =- \mu
^2 \partial _\omega \|   \phi \| _2^2
\endaligned $$ by $H_\omega
\partial _\omega\Phi =-\sigma _3 \Phi$. So $N(\sigma _3H_\omega )\le 1$, which is
incompatible with $N(\sigma _3H_\omega )\ge 2.$ To conclude the
proof of Theorem 1.1 we need to prove
   (2.3):

   \proclaim{Proposition 2.4} If $H_\omega
$ satisfies the conditions in Definition 2.3 for any non zero $u\in
H^1\cap N^\perp _g (H_\omega ^\ast )$ we have $\langle \sigma
_3H_\omega u, u \rangle \ge 0$.
\endproclaim
Under our hypotheses we have the decomposition
$$ \align  &   N_g^\perp (H^\ast _{\omega })=
\sum  _{\lambda \in \sigma _p\backslash \{ 0\}  }\ker (H _{\omega }
- \lambda   )  \oplus L_c^2(H_{\omega }) \tag 1 \\& \text{ where }
L_c^2(H _\omega )=  \{   N_g(H ^\ast _\omega ) \oplus \sum _{\lambda
\in \sigma _p\backslash \{ 0\} }\ker (H ^\ast _\omega -\lambda
 )  \} ^{\perp}.\endalign $$
We have $\langle\sigma _3 H_\omega u ,  v \rangle =0$  for $u$ and
$v$ in different terms in (1). By hypothesis, $\langle\sigma _3
H_\omega \cdot , \cdot \rangle  $ is a positive quadratic form in
each $\ker (H _{\omega } - \lambda   )  .$ So Proposition 2.4 is a
consequence of
$$\langle \sigma _3H_\omega u, u \rangle \ge 0 \text{ for any $u\in
H^1\cap L_c^2(H_{\omega })$.} \tag 2.4 $$ (2.4) is a general fact.
In this section we will consider some special cases for $H_\omega$
and in \S 5 we will complete the proof of (2.4). First of all we
remind the following definition: \proclaim{Definition 2.5} $\omega $
is a resonance if there is a distribution $F$ such that $H_\omega
F=\omega F$ such that $F\in L^{2, -s}$ for any $s>1/2$ but $F\not
\in L^2$.
\endproclaim
If $\omega $ is neither a resonance nor an eigenvalue of $H_\omega$,
which is a generic condition, then by Theorem 2.11 \cite{CPV} in
(2.4) we have $\langle \sigma _3H_\omega u, u \rangle
> 0$ for $u\neq 0$.
We consider now the case when $\omega $ is a    resonance or  an
eigenvalue.  Using the  terminology in Jensen and Kato \cite{JK}  we
can distinguish between $\omega$ being  exceptional point of first
kind   (when $\omega $ is a resonance but not an eigenvalue), second
kind (when $\omega$ is an eigenvalue but not a resonance) and third
kind ($\omega$ both resonance and eigenvalue). In this section we
consider two special cases. The proof is then completed in \S 5.

\proclaim{Lemma  2.6}  Let
 $f,g \in C^\infty _0(\Bbb R^3,\Bbb R)$ and set $U_1(x)=\sigma _3
 f   +i g
 \sigma _2.$ Suppose that in the space $\Cal V$ formed by eigenfunctions and
 resonant functions at $\omega$  the quadratic form $$\langle \sigma _3U_1 \cdot , \cdot  \rangle
 \text{ is strictly positive.} \tag 2.5$$
Suppose furthermore that $\omega $ is either exceptional    of first
kind or of second type. If $\omega $ is      of second type assume
furthermore that $\dim \Cal V=1$. Then for any $\varepsilon
>0$ sufficiently close to 0 we have:

{\item {(1)}}  the point $\omega$ is neither a resonance nor an
eigenvalue for $H_{\omega, \varepsilon }:=H_\omega +\varepsilon
U_1$.

{\item {(2)}}  $H_{\omega, \varepsilon }$ does not have eigenvalues
close to $\omega$.
\endproclaim
Lemma 2.6 is proved \cite{CuP}. Notice that  $\dim \Cal V<\infty$
and probably the statement holds always without the restriction
$\dim \Cal V=1$. Nonetheless, in \S 5 we give a different proof of
the remaining cases (2.4). Assuming the conclusions of Lemma 2.6,
which are valid for $\dim \Cal V=1$, by $\dim \Cal V<\infty$ whe
know that there are $U_1$ as above satisfying (2.5). Let now $\gamma
$ be a fixed and small counterclockwise circle with center the
origin in $\Bbb C$.
$$\aligned  \text{Set }\quad & P_{\varepsilon}:=-\frac{1}{2\pi i}\int _{\gamma}R_{H_{\omega,
\varepsilon }}(z) dz =P ^{(0)} -\varepsilon P^{(1)}+\varepsilon
^2P(\varepsilon  )\\& P^{(0)} \text{ is the projection on
$N_g(H_\omega )$}
\\& P^{(1)}=
\frac{ 1}{2\pi i}\int _{\gamma}R_{H_{\omega  }}(z)U_1 R_{H_{\omega
}}(z)dz \\& P(\varepsilon  ) = \frac{ 1}{2\pi i}\int
_{\gamma}R_{H_{\omega  }}(z)U_1 R_{H_{\omega, \varepsilon }}(z)U_1
R_{H_{\omega }}(z) dz . \endaligned
$$
Then $\| P^{(1)}:L^2\to H^1\|+ \| P(\varepsilon  ):L^2\to H^1\| <C$
for a fixed constant if $ \varepsilon \in [0,\varepsilon _0]$  for
some $\varepsilon _0>0$ small enough. So, if we consider $u\in
H^1\cap L_c^2(H_{\omega })$ and we split $ u=u_1+u_2$ with
$u_1=P_{\varepsilon}u$, we have $u_1= O(\varepsilon )$ and
$$\aligned & \langle \sigma _3H_\omega u, u \rangle =\langle \sigma _3H_{\omega,
\varepsilon } u, u \rangle -\varepsilon \langle \sigma _3U_1u,
u\rangle = \langle \sigma _3H_{\omega, \varepsilon } u_1, u_1
\rangle + \\& + \langle \sigma _3H_{\omega, \varepsilon } u_2, u_2
\rangle  -\varepsilon \langle \sigma _3U_1u, u\rangle > \langle
\sigma _3H_{\omega, \varepsilon } u_1, u_1 \rangle   -\varepsilon
\langle \sigma _3U_1u, u\rangle =O(\varepsilon )
 \endaligned
$$
and so $\langle \sigma _3H_\omega u, u \rangle \ge 0.$

We consider the case $\dim \Cal V>1$ in \S 5.

\head \S 3  Orbital instability of excited states: the case when the
internal modes are close to the continuous spectrum
\endhead

  We will assume the following hypotheses.

{\item {(H4)}} $\beta (t)\in C^3(\Bbb R, \Bbb R)$.

{\item {(H5)}} The operators $  L_{+,\omega} =-\Delta + \omega
-\beta (\phi ^2_\omega ) -2\phi ^2 _\omega \beta ^\prime (\phi
^2_\omega  ) $ are such that $\ker L_{+,\omega}\cap H^1_{r }(\Bbb
R^3)=0 $ with $H^1_{r }$  introduced in (H1).

{\item {(H6)}}  In (1.1) the initial data are $u_0\in H^k(\Bbb R^3)$
for $k=2$ and satisfy  $u_0(x)=u_0(-x)$. Let $ H^k_e(\Bbb R^3)$ the
space of such functions.

{\item {(H7)}}  $\frac d {d\omega } \| \phi _ {\omega
}\|^2_{L^2(\Bbb R ^3 )}\neq 0 $  for $\omega\in\Cal{O}$ and $N_g(H
_\omega )\cap H^k_e$ be is spanned by $\{ \sigma _3\Phi _\omega ,
\partial _\omega \Phi _\omega \} $;

{\item {(H8)}} Let $H_\omega$ be the linearized operator around
$e^{it\omega}\phi_\omega$, see Section 2. Then $H_\omega$ has a
certain number of
  simple positive eigenvalues with $0< \lambda _j(\omega )<
\omega < 2\lambda _j(\omega )$. $H_\omega $ does not have other
eigenvalues and $\pm \omega $ are not resonances.

{\item {(H9)}} For   multi indexes  $m=(m_1,m_2,...)$ and
$n=(n_1,...)$, setting $\lambda (\omega )=(\lambda _1(\omega ),...)$
and $(m-n)\cdot \lambda =\sum (m_{j}-n_{j})  \lambda _j  $, we have
the following two non resonance hypotheses: {\item {(i)}}
$(m-n)\cdot \lambda (\omega )=0$ implies $m=n$ if $|m|\le 3$ and
$|n|\le 3$;  {\item {(ii)}} $(m-n)\cdot \lambda (\omega )\neq \omega
$ for all $(m,n)$ with $|m|+|n|\le 3$.

{\item {(H10)}} We assume the non degeneracy  Hypothesis 3.7.

The key hypotheses are (H8), where the condition $\lambda _j(\omega
)>\omega /2$ is a quantitative description of what it means for the
eigenvalues to be close to the  continuous spectrum, and (H10),
which is valid generically.

 A standing wave is orbitally unstable if it is
not orbitally stable. Recall the following definition:

\proclaim{Definition 3.1}  A standing wave $e^{i\omega t} \phi _
{\omega } (x)$ is orbitally stable if for any $ \epsilon >0$ there
is a $\delta (\epsilon )>0$ such that   for any  $\| u(0,x) - \phi _
{\omega } \| _{H^1_{x}}<\delta  (\epsilon )$ the corresponding
solution $u(t,x)$ is globally defined and for any $t$ we have
$$\inf _{\gamma \in \Bbb R \, \& \, x_0 \, \in \Bbb R^n }
  \|  u(t,x) -e^{i \gamma }\phi _ {\omega } (x-x_0) \| _{H^1 _x}
  <   \epsilon .
$$\endproclaim
{\it Remark.} In the setting of even solutions of (1.1) or of
solutions of the non translation invariant (1.3), we need to  pick
$x_0=0$ in the above definition.

In this section we will prove:

\proclaim{Theorem 3.2 } Under hypotheses (H1-10) the excited states
$e^{i\omega t}\phi _{\omega }(x)$ are orbitally unstable.

\endproclaim
{\it Remark.} It is easy to manufacture examples in the spirit of
\cite{SW1-3,TY1-3,T} by considering short range Schr\"odinger
operators $-\Delta +q(x)$ which admit a certain number of
eigenvalues $-E_0<-E_1<...<0$. For example, if there are only two
simple eigenvalues   with $E_1/2<E_0-E_1<E_1$, then if $q(x)$ is
generic  the  excited states originating from $E_1$ are unstable.

The proof of Theorem 3.2  covers the reminder of \S 3. We assume by
absurd that the excited state   $e^{i\omega _0 t} \phi _ {\omega _0}
(x)$ is orbitally stable. We pick an arbitrarily small $\epsilon >0$
and we consider the associated $\delta =\delta (\epsilon )< \epsilon
$.
 Then  the  representation (2.1) is valid for all $t\in \Bbb R$
 with:     $r(t,x)\in
C(\Bbb R, H^2 (\Bbb R^3))$, $R(t)\in N_g^\perp (H^\ast _{\omega
(t)}) $ for all $t$, $|\omega _0-\omega (0)|<C\delta $ for fixed
$C>0$ and $|\omega _0-\omega (t)|<C\epsilon $ for all $t$. Recall
that there is a real valued function $F(|u|^2)$ with $\beta (|u|^2)
u=\partial _{\overline{u}}( F(|u|^2))$ and $F(0)=0.$ For $u$ given
by (2.1), we have $e^{-i\vartheta }\beta (|u|^2) u=\partial
_{\overline{r}}F(|\phi _\omega +r|^2)$. Recall
$$ \aligned & e^{-i\vartheta }\beta (|u|^2) u=  \beta (\phi _{\omega  } ^2) \phi _{\omega  }
+\left (\beta   ( \phi _{\omega  } ^2 )    +\beta ^\prime ( \phi
_{\omega (t)} ^2 )\phi _{\omega  } ^2\right ) r
 +\beta ^\prime ( \phi _{\omega  } ^2 )
\phi _{\omega  } ^2  \overline{  r }- n(r,\overline{r} ).\endaligned
$$
We have $\left (\beta   ( \phi _{\omega  } ^2 )    +\beta ^\prime (
\phi _{\omega (t)} ^2 )\phi _{\omega  } ^2\right ) r
 +\beta ^\prime ( \phi _{\omega  } ^2 )
\phi _{\omega  } ^2  \overline{  r }=$
$$\aligned &    \frac 1{2}\partial
_{\overline{r}} \left (   \left ( \beta   ( \phi _{\omega  } ^2 )
+\beta ^\prime ( \phi _{\omega  } ^2 )\phi _{\omega  } ^2\right )
{|r|^2}
 +\beta ^\prime ( \phi _{\omega  } ^2 )
\phi _{\omega  } ^2  \left ( \overline{  r }^2 +r^2 \right ) \right
)
\endaligned $$
Then   $- n(r,\overline{r} )=\partial _{\overline{r}} G(R )$ with
$G(R )$ real valued, $G(0)=0$. For $^t  G'(R)=(
G_{r},G_{\overline{r}})$
$$i  R _t =H _{\omega (t)}   R +\dot \gamma (t) \sigma _3R+\sigma _3 \sigma _1 G'(R)
 +\dot \gamma (t) \sigma _3\Phi _\omega
-i \dot \omega (t) \partial _\omega \Phi _\omega . \tag 3.1
$$ The condition $R(t)\in N_g^\perp (H^\ast _{\omega (t)}) $ yields
the modulation equations
$$\aligned & i\dot \omega \langle \Phi  , \partial
_\omega \Phi
 \rangle
=\langle \sigma _3 \dot \gamma    R+N(R)+i\dot \omega \partial
_\omega P _{N_g(H_\omega )}R, \Phi \rangle
\\& \dot \gamma \langle \Phi  , \partial _\omega \Phi
 \rangle
=-\langle \text{same as above}, \sigma _3 \partial _\omega \Phi
\rangle .
\endaligned  \tag 3.2$$
These can be used to express $i\dot \omega =i\dot \omega (\omega ,
R)$, $ \dot \gamma = \dot \gamma (\omega , R)$.

\proclaim{Lemma 3.3} We can write with  smooth
 functions in
 $r$ and $\overline{r}\in H^{2}(\Bbb R ^3)$
 $$\aligned & i\dot \omega =i\dot \omega (\omega ,r,\overline{r}) =
 \nu (r,\overline{r}) -\nu (\overline{r}, {r})
  \text{ with $\overline{\nu (r,\overline{r})}=\nu (\overline{r}, {r})$} \\&
 \dot \gamma =\dot \gamma (\omega ,r,\overline{r})=
 \mu (r,\overline{r}) +\mu (\overline{r}, {r})
 \text{ with $\overline{\mu (r,\overline{r})}=\mu (\overline{r},
 {r})$}.\endaligned \tag 1
 $$

\endproclaim
{\it Proof.} Let $ P_{N_g(H_\omega )}$ be the projection onto
$N_g(H_\omega )$ in $L^2=N_g(H_\omega )\oplus N_g^{\perp }(H_\omega
^\ast)$. We apply   $ P_{N_g(H_\omega )}$ to (2.2) obtaining
$$i  P_{N_g(H_\omega )}R _t - \dot \gamma (t) P_{N_g(H_\omega )}\sigma _3   R
-\sigma _3 \dot \gamma (t)\Phi _\omega+i \dot \omega (t)\partial
_\omega \Phi _\omega+P_{N_g(H_\omega )}N(R )   .$$ Set $q(\omega)=\|
\phi _\omega \| ^2_2$ and $q'(\omega) =dq(\omega)/d\omega$. Then we
have
$$ P_{N_g(H_\omega )} =\sigma _3\Phi _\omega \langle \quad , \sigma
_3\partial _\omega \Phi _\omega \rangle /q(\omega )+\partial _\omega
\Phi _\omega \langle \quad ,   \Phi _\omega \rangle /q(\omega ) .$$
By $P_{N_g(H_\omega )}R=0 $, which implies $ P_{N_g(H_\omega )}R
_t=-\dot \omega  \partial _\omega P _{ N_g (H _\omega )}  R$, we get

$$\aligned & \left (  q'(\omega) +\left [    \matrix -\langle  \partial _\omega P _{ N_g
(H   _\omega )}  R , \Phi _\omega \rangle
&    \langle \sigma _3 R , \Phi _\omega \rangle  \\
- \langle  \partial _\omega P _{ N_g  (H   _\omega )}  R, \sigma _3
\partial _\omega \Phi _\omega \rangle  & \langle   R , \partial
_\omega \Phi _\omega \rangle
\endmatrix \right ] \right )     \left [ \matrix  i\dot \omega \\ -\dot
\gamma
\endmatrix  \right ] = \left [ \matrix  \langle N(R), \Phi _\omega \rangle  \\
\langle N(R), \sigma _3 \partial _\omega \Phi _\omega \rangle
\endmatrix \right ] .
\endaligned  $$

By an elementary computation we have
$$ \aligned & \langle  \partial _\omega P _{ N_g
(H   _\omega )}  R , \Phi _\omega \rangle =\left \langle
r+\overline{r}, \frac{\langle \partial ^2_\omega \Phi ,\Phi \rangle
}{q'} \phi _\omega + q'
\partial _\omega \frac{\phi _\omega }{q'(\omega )} \right \rangle \\&\langle  \partial _\omega P _{ N_g
(H   _\omega )}  R , \sigma _3 \partial _\omega \Phi _\omega \rangle
=\left \langle r-\overline{r}, \frac{\| \partial  _\omega \Phi \|
_2^2 }{q'} \partial _\omega \phi _\omega + q'
\partial _\omega \frac{\partial _\omega \phi _\omega }{q'(\omega )}\right  \rangle
\endaligned $$
and so, for some real valued exponentially decreasing functions
$\alpha (\omega ,x)$ and $\beta (\omega ,x)$, we have the following,
which yields (1):

$$\aligned & \left (   q'(\omega )   +\left [    \matrix  \langle   r+\overline{r}, \alpha (\omega ) \rangle
&   \langle r-\overline{r} , \phi _\omega \rangle  \\
 \langle r-\overline{r}, \beta (\omega ) \rangle  & \langle
r+\overline{r} , \partial _\omega \phi _\omega \rangle
\endmatrix \right ] \right )     \left [ \matrix  i\dot \omega \\ -\dot
\gamma
\endmatrix  \right ]  = \left [ \matrix   \langle n(r,\overline{r}) -n( \overline{r},r), \phi _\omega \rangle  \\
  \langle n(r,\overline{r})+n( \overline{r},r),   \partial _\omega \phi _\omega
\rangle
\endmatrix \right ] .
\endaligned  $$
The regularity of $\mu ( r, \overline{r}) $ and $\nu ( r,
\overline{r}) $ follows by the smoothness of $n(z,\overline{z})$ as
a function in $z\in \Bbb C$, and by the fact that $H^2(\Bbb R^3) $
is an algebra. This completes Lemma 3.3.
\bigskip

For each $j$ we consider a generator $\xi _j\in \ker (H_\omega
-\lambda _j)$ such that $\langle \xi _j, \sigma _3 \xi _j \rangle =s
_j$ with
 $s _j=1$ (resp. $s _j=-1$) if $\lambda _j$ has positive
 (resp. negative) signature. Since $e^{i\omega t} \phi _ {\omega } $ is an excited state, by
Theorem 1.1 at least for one $j$ we have $s _j=-1$, so in particular
we can assume $s _1=-1$. Indeed under hypothesis (H8), if $\xi$ is a
generator of $\ker (H_\omega -\lambda _j)$ for any $j$, then
$\langle \xi , \sigma _3\xi \rangle \neq 0$. We expand $R(t)\in
N_g^\perp (H^\ast _{\omega (t)}) $ into

$$    R (t) =(z\cdot \xi + \bar z\cdot \sigma _1 \xi ) + f(t)   \in
\big [ \sum _{j,\pm  }\ker (H  _{\omega (t)}\mp \lambda _j(\omega
(t)))\big ] \oplus L_c^2(H_{\omega (t)})  . \tag 3.3$$
Correspondingly we   express  (3.1) as

$$\aligned   i\dot z _j\xi _j-\lambda _j(\omega ) z_j\xi _j&= P_{\ker (H  _\omega -\lambda _j)}
 ( \dot \gamma (\omega , R) \sigma _3R+\sigma _3 \sigma _1 G'(R)\\&
 -iz_j\dot \omega (\omega , R)
 \partial _\omega \xi _j +i\dot \omega (\omega , R) \partial _\omega P_{\ker (H  _\omega -\lambda _j)}R)
  \\  iP_c(H_\omega ) \dot f -H_\omega f &= P_c(H_\omega )
(  \dot \gamma (\omega , R) \sigma _3R+\sigma _3 \sigma _1
G'(R)+i\dot \omega (\omega , R) \partial _\omega P_{c} (H  _\omega
)R).\endaligned \tag 3.4$$ We   use the multi index notation
$z^m=\prod _j z^{m_j}_j$. We  consider the expansion
$$\aligned & \sigma _3 \sigma _1 G'(R)=\sum _{   |m+n|=2}^{ 3} R_{m,n}(\omega ) z^m  \bar z^n+
\sum _{ |m + n|=1}  z^m  \bar z^n A_{m,n}(\omega ) f+ O(f^2)+\cdots
\endaligned  \tag  3.5$$ with $R_{m,n}(\omega  ,x) $ and
$A_{m,n}(\omega ,x ) $ real vectors  and matrices exponentially
decreasing in $x$. We have
$$\aligned & A_{m,n}(\omega )= \frac{\sigma _3 \sigma _1}{m! n!}\partial _z^m
 \partial _{\overline{z}} ^n\partial _f G'(0)
  \, , \quad  R_{m,n}(\omega )= \frac{\sigma _3
\sigma _1}{m! n!}\partial _z^m  \partial _{\overline{z}} ^nG'(0)
 .\endaligned$$
Notice $A_{m,n}(\omega )=-\sigma _1 A_{n,m}(\omega )\sigma _1 $ and
$\sigma _1R_{m,n}(\omega )=-R_{n,m}(\omega )$. Indeed  by definition
$\sigma _1 R=\overline{R}$, from which we get $\sigma _1
f=\overline{f}$ and $\sigma _1 G'(R)=\overline{G'(R)}$. We also have
$\sigma _1H_\omega=-H_\omega \sigma _1$. Then, taking complex
conjugate of (3.1) and applying to the resulting equation $\sigma
_1$, we get
$$\aligned & - i\dot R+H_\omega R = ..\sigma _1R_{m,n}(\omega )z^{n}
\overline{z}^{m}+.. \overline{z}^mz^n\sigma _1A_{m,n}(\omega )\sigma
_1 f+.. \\& =..- R_{n,m}(\omega )z^{n} \overline{z}^{m}+..
-\overline{z}^mz^n A_{n,m}(\omega )  f+.. \endaligned$$ which yields
$A_{m,n}(\omega )=-\sigma _1 A_{n,m}(\omega )\sigma _1 $ and $\sigma
_1R_{m,n}(\omega )=-R_{n,m}(\omega )$. We set $\delta _j=(\delta
_{j1},\delta _{j2},...) $
 with $\delta _{jk}$ the Kronecker delta. We
have
$$ \aligned &  A_{\delta _\ell ,0}(\omega ) =
\sigma _3 \sigma _1 \partial _{z_\ell}
 \partial _f G'(0)  = \sigma _3 \sigma _1\partial
_{z_\ell }
 \partial _f G'(0)  =  \sigma _3 \sigma _1 G ^{(3)}(0)
    ( \quad ,\xi _\ell , P_c(H  _\omega  ))  \endaligned $$
where $G ^{(3)}(0)$ is written as a symmetric trilinear form and
where one of the vectors of the triple is $\xi _\ell$. We have

$$ \aligned & P_c(H  _\omega  )R_{ \delta _j +\delta _\ell ,0}(\omega ) =
\frac{\sigma _3 \sigma _1}{(\delta _j +\delta _\ell )!}P_c(H
_\omega  ^\ast )\partial _{z_j}\partial _{z_\ell }
   G'(0)  \\& = \frac{\sigma _3 \sigma _1}{(\delta _j +\delta _\ell )!}\partial
_{z_j}\partial _{z_\ell }
   G'(0) \circ P_c(H  _\omega    )=  \frac{\sigma _3 \sigma _1}{(\delta _j +\delta _\ell )!} G ^{(3)}(0)
 (\xi _j,\xi _\ell  ,  P_c(H  _\omega  ))  .\endaligned $$
For later use we record:

\proclaim{Lemma 3.4} For $\lambda _j=\lambda _j(\omega )$,  $\xi
_j=\xi _j(\omega )$,  $A_{0, \delta _\ell  }  =  A_{0, \delta _\ell
} (\omega )$ and  $R_{ \delta _j +\delta _\ell ,0}  =  R_{ \delta _j
+\delta _\ell ,0} (\omega )$ we have
$$\aligned &  \langle A_{0, \delta _\ell  } R_{H_\omega }^+(\lambda  _j +\lambda
_\ell )R_{ \delta _j +\delta _\ell ,0}  , \sigma _3 \xi _j \rangle =
\frac{1}{(\delta _j +\delta _\ell )!}\times \\&   \langle
 R_{H_\omega   }^+(\lambda _j +\lambda _\ell )
  \sigma _3\sigma _1G ^{(3)}(0)
 (\xi _j,\xi _\ell  ,  P_c(H  _\omega  ))
  ,  \sigma _3 \sigma _3\sigma _1G ^{(3)}(0)
 (\xi _j,\xi _\ell  ,  P_c(H  _\omega  )) \rangle .
\endaligned $$
Taking the imaginary part of the above formula, we have

$$\aligned &  \Im \langle A_{0, \delta _\ell  } R_{H_\omega }^+(\lambda  _j +\lambda
_\ell )R_{ \delta _j +\delta _\ell ,0} (\omega ), \sigma _3 \xi
_j(\omega ) \rangle =    \frac{\pi }{(\delta _j +\delta _\ell
)!}\times \\& \langle
 \delta ( H_\omega   -\lambda _j - \lambda _\ell )
  \sigma _3\sigma _1G ^{(3)}(0)
 (\xi _j,\xi _\ell  ,  P_c(H  _\omega  ))
  ,  \sigma _3 \sigma _3\sigma _1G ^{(3)}(0)
 (\xi _j,\xi _\ell  ,  P_c(H  _\omega  )) \rangle \\& \ge 0 \text{ for any $G ^{(3)}(0)
 (\xi _j,\xi _\ell  ,  P_c(H  _\omega  )).$}
\endaligned $$

\endproclaim
{\it Proof.} We recall  that  $G ^{(3)}(0)
 (\xi _j,\xi _\ell  ,  P_c(H  _\omega  ))\in L^2$ is defined by the
 equality
 $\langle g, G ^{(3)}(0)
 (\xi _j,\xi _\ell  ,  P_c(H  _\omega  ))\rangle =\frac{d}{dt}G ^{(2)}(tP_c(H  _\omega  )g)_{t=0}
 (\xi _j,\xi _\ell   ).$
We assume the first formula and we set $F=\sigma _3\sigma _1G
^{(3)}(0)
 (\xi _j,\xi _\ell  ,  P_c(H  _\omega  ))$. Recall
there are isomorphisms $W(H_\omega   ):L^p \to L^p_c(H_\omega )  $
and $Z(H_\omega   )$ its inverse, such that $P _{c}(H_\omega
)H_\omega   =W(H_\omega   ) \sigma _3(-\Delta +\omega )Z(H_\omega
  ) $, \cite{Cu2,CPV}. Furthermore, from the definitions one
gets   $ W^\ast (H_\omega   ) \sigma _3= \sigma _3Z (H_\omega )$,
see \cite{Cu2,CPV}. For $P _{c}(H_\omega   )F=W(H_\omega
)\widetilde{F}$ we have
 $$\aligned &  \langle
 R_{H_\omega   }^+(\lambda _j +\lambda _\ell )
  F
  ,  \sigma _3 F \rangle =\langle
 R_{H_\omega   }^+(\lambda _j +\lambda _\ell )
 W(H_\omega   )\widetilde{F}
  ,  \sigma _3 W(H_\omega   )\widetilde{F} \rangle  =\\  =& \langle
 R_{\sigma
_3(-\Delta +\omega )   }^+(\lambda _j +\lambda _\ell )
 \widetilde{F}
  , W^\ast (H_\omega   ) \sigma _3 W(H_\omega   )\widetilde{F}
  \rangle = \\  = &\langle
 R_{\sigma
_3(-\Delta +\omega )   }^+(\lambda _j +\lambda _\ell )
 \widetilde{F}
  ,  \sigma _3Z (H_\omega   ) W(H_\omega   )\widetilde{F}
  \rangle = \langle
 R_{\sigma
_3(-\Delta +\omega )   }^+(\lambda _j +\lambda _\ell )
 \widetilde{F}
  ,  \sigma _3 \widetilde{F}
  \rangle
 \endaligned $$
Let now $\widetilde{F}_1$ and $\widetilde{F}_2$ be the two
components of the vector $ \widetilde{F} $. Then
$$ \aligned &\Im \langle
 R_{\sigma
_3(-\Delta +\omega )   }^+(\lambda _j +\lambda _\ell )
 \widetilde{F}
  ,  \sigma _3 \widetilde{F}
  \rangle =\pi \langle
 \delta ( \sigma
_3(-\Delta +\omega ) - \lambda _j -\lambda _\ell )
 \widetilde{F}
  ,  \sigma _3 \widetilde{F}
  \rangle \\  =&   \pi \langle
 \delta (  -\Delta +\omega   - \lambda _j -\lambda _\ell )
 \widetilde{F}_1
  ,    \widetilde{F}_1
  \rangle - \pi \langle
 \delta ( \sigma
 \Delta - \omega   - \lambda _j -\lambda _\ell )
 \widetilde{F}_2
  ,    \widetilde{F} _2
  \rangle . \endaligned
$$
We have $\langle
 \delta (
 \Delta - \omega   - \lambda _j -\lambda _\ell )
 \widetilde{F}_2
  ,    \widetilde{F} _2
  \rangle =0$ for any $ \widetilde{F} _2 $.
Notice that in our hypothesis we have $\lambda _j +\lambda _\ell
>\omega . $ Setting $\rho _0=\sqrt{\lambda _j +\lambda _\ell -\omega  }$
we have for any $ \widetilde{F} _1 \in H^{s}$ for $s>1/2$
$$\langle
 \delta (  -\Delta +\omega   - \lambda _j -\lambda _\ell )
 \widetilde{F}_1
  ,    \widetilde{F}_1
  \rangle =  \frac{1}{2\rho _0}\int _{|\eta |=\rho _0} |\widehat{\widetilde{F}_1}(\eta )|^2
  d\sigma (\eta )\ge 0.\tag 3.6$$
We prove now the first formula in the statement. We have
$$\aligned &  \langle A_{0, \delta _\ell  } R_{H_\omega }^+(\lambda  _j +\lambda
_\ell )R_{ \delta _j +\delta _\ell ,0}, \sigma _3 \xi _j \rangle = -
\langle \sigma _1A_{  \delta _\ell ,0  } \sigma _1R_{H_\omega
}^+(\lambda _j +\lambda _\ell )R_{ \delta _j +\delta _\ell ,0} ,
\sigma _3 \xi _j \rangle =\\& = - \langle \sigma _1\sigma _3 \sigma
_1 \partial _{z_\ell }\partial _{f }
   G'(0) \sigma _1R_{H_\omega }^+(\lambda _j
+\lambda _\ell )R_{ \delta _j +\delta _\ell ,0} , \sigma _3 \xi _j
\rangle =\\& =   \langle  \partial _{z_\ell }\partial _{f }
   G'(0) \sigma _1R_{H_\omega }^+(\lambda _j
+\lambda _\ell )R_{ \delta _j +\delta _\ell ,0} ,   \xi _j \rangle
\\& =\langle    \sigma _1R_{H_\omega }^+(\lambda _j +\lambda _\ell )R_{
\delta _j +\delta _\ell ,0} ,  G ^{(3)}(0)
 (\xi _j,\xi _\ell  ,  P_c(H  _\omega  )) \rangle =\\& =\frac{1}{(\delta _j +\delta _\ell )!}\langle    \sigma _1
 R_{H_\omega }^+(\lambda _j +\lambda _\ell )
 \sigma _3 \sigma _1G ^{(3)}(0)
 (\xi _j,\xi _\ell  ,  P_c(H  _\omega  ))
  ,  G ^{(3)}(0)
 (\xi _j,\xi _\ell  ,  P_c(H  _\omega  )) \rangle  \\& =\frac{1}{(\delta _j +\delta _\ell )!}\langle
 R_{H_\omega   }^+(\lambda _j +\lambda _\ell )h _{j,\ell },
  \sigma _3h _{j,\ell } \rangle   \text{ where $h _{j,\ell }=\sigma _3\sigma _1G ^{(3)}(0)
 (\xi _j,\xi _\ell  ,  P_c(H  _\omega  ))$.}\endaligned
$$
 This concludes the proof of Lemma 3.4.

\bigskip

For  $O_{loc}(z^n)=\sum _\ell O_{loc}(|z_\ell ^{n}|)$, (3.4) can be
expressed as

$$\aligned &if_t=\left ( H _{ \omega (t)}+P_c(H_\omega )\sigma _3 \dot \gamma \right )f  +
\sum _{  |m+n|=2} z^m \bar z^nP_c(H_\omega )R_{m,n}(\omega )
\\&  + \sum _{  |m + n|=1} z^m \bar z^n P_c(H_\omega )A_{m,n}(\omega ) f+
O(f^2)+O_{loc}(z^3),
\endaligned \tag  3.7
$$
 and
$$\aligned i\dot z _j\xi _j-\lambda _j(\omega ) z_j\xi _j&= P_{\ker (H  _\omega -\lambda _j)}
 ( \sum _{  |m+n|=2} z^m \bar
z^nR_{m,n}(\omega )+\sum _{  |m+n|=3} z^m \bar
z^nR_{m,n}^{(1)}(\omega )\\&
 +\sum _{  |m + n|=1} z^m \bar z^n A_{m,n}(\omega ) f+O(f^2)+O_{loc}(z^4)
\endaligned \tag  3.7
$$
where in (3.7) the coefficients $A_{m,n}(\omega )$ and
$R_{m,n}(\omega )$ are those of the expansion of $\sigma _3 \sigma
_1 G'(R)$ in (3.5) and where $R_{m,n}^{(1)}(\omega )$ are real
  and exponentially decreasing vectors. We have:

  \proclaim{Lemma 3.5} Let $0<\delta < \epsilon $ be as in the Definition
  3.1 of linear stability. For any $C_1>0$
there are a $\varepsilon (C_1)>0$ and a $C(C_1)$ such
  that if, for $0<\varepsilon <\varepsilon (C_1)$, we have
  $\| z_j\| _{L^4(0,T)}^2\le C_1\varepsilon $ for all $j$, then
for all admissible pairs $(p,q)$ we have for a fixed $c_0$
$$\| f\| _{L^p_t((0,T), W^{1,q}_x)}<c_0\delta +C(C_1)\varepsilon .\tag 3.8$$
\endproclaim
{\it Proof.} In this proof $P_c(\omega )=P_c(H _{\omega} )$. We
split $P_c(\omega )= P_+(\omega )+P_-(\omega )$, with  $P_\pm
(\omega )$ the spectral projections in $\Bbb R_\pm \cap \sigma
_c(H_\omega )$, see \cite{Cu3}. By orbital stability we can  fix
$\omega _0$ such that $|\omega (t)-\omega _0|=O(\epsilon )$ for all
$t$. We write following \cite{BP}
$$\aligned &if_t=\left \{ H_ {\omega _0}   + (\dot \gamma +\omega
-\omega _0) (P_+(\omega _0)-P_-(\omega _0))\right \} P_c(\omega _0)f
\\& +O_{loc}(\epsilon f)+
O(f^2)+O_{loc}(z^2)\\& \text{where } O_{loc}(\epsilon f)=
 (\dot \gamma +\omega -\omega _0) \left (P_c(\omega _0)\sigma _3-  (P_+(\omega _0)-P_-(\omega _0))  \right ) f
 \\&  +
   \left (  V  (\omega ) -
 V  (\omega _0) \right )   f
  +(\dot \gamma +\omega -\omega _0) \left ( P_c(\omega  )-  P_c(\omega _0)\right )\sigma _3f
  .
\endaligned \tag 3.9
$$
To justify the notation $O_{loc}(\epsilon f)$ we notice that $\omega
-\omega _0=O(\epsilon )$ and the following fact: $\forall$ $p \in
[1, 2]$ $q \in [2, \infty )$ with $ c_{p,q} (\omega )$ upper
semicontinuous in $\omega$,    \cite{Cu3,Cu1},

 $$ \|   P_c (\omega
) \sigma _3- (P_+ (\omega )-P_- (\omega ) ) : L^q\to L^p  \|  \le
c_{p,q} (\omega )<\infty .$$ By \cite{Cu2,CPV} $ P_c (\omega
)e^{-itH_\omega }$ satisfies for any fixed $\omega \in \Cal O$ the
Stricharz estimates, i.e. there is a $C(\omega ,k)$ upper
semicontinuous in $\omega$ such that for all admissible pairs
$(p,q)$ and $(a,b)$ and
  we have

$$ \aligned & \| P_c (\omega )e^{-itH_\omega }\varphi (x) \|
_{L^p_tW^{k,q}_x}\le C(\omega ) \|  \varphi (x) \| _{ H^{k }_x} \\&
\| P_c (\omega )e^{-itH_\omega }\psi  \| _{L^p_tW^{k,q}_x} \le
C(\omega ) \|
 \psi (t,x) \| _{L^{a'}_tW^{k,b'}_x}.\endaligned \tag 3.10 $$
We have $\| f \| _{L^p_tW^{k,q}_x}\approx  \| P_c(\omega _0)f \|
_{L^p_tW^{k,q}_x}$ by

$$\aligned & \| f \| _{L^p_tW^{k,q}_x}= \| P_c(\omega _0)f
\| _{L^p_tW^{k,q}_x} +\| \left ( P_c(\omega  )-  P_c(\omega
_0)\right )f \| _{L^p_tW^{k,q}_x}\\& = \| P_c(\omega _0)f \|
_{L^p_tW^{k,q}_x} + O(\epsilon)\|   f \|
_{L^p_tW^{k,q}_x}\endaligned $$ For $  U_\pm (t,t')= e^{-i(t-t')H
_{\omega_0}} e^{\pm i\int _{t'}^t d\tau (\dot \gamma (\tau )+\omega
(\tau )-\omega _0) }P _{\pm}(\omega _0)$, we have
$$\aligned & P _{\pm}(\omega _0)f(t)=  U_\pm (t,0)f(0)+\int _0^tU_\pm (t,t')
(O_{loc}(\epsilon f)+ O(f^2)+O_{loc}(z^2)) dt'.
\endaligned $$
Since $\| f(0) \| _{H^1}<\delta $, there is a $T_1\in (0,T]$ such
that (3.8) is true in $(0,T_1)$. Using the Stricharz estimates
(3.10), in particular the "endpoint Stricharz estimate, in $(0,T_1)$
we have
$$\aligned & \| f \| _{L^p_tW^{k,q}_x}\approx  \| P_c(\omega _0)f \|
_{L^p_tW^{k,q}_x}\\& \le C(\omega )\delta +  C(\omega )\epsilon \| f
\| _{L^2_tW^{k,6}_x}+ C(\omega ) \epsilon O(\delta ^2 +\varepsilon
^2)+C(\omega )C_1\varepsilon .
\endaligned $$
Since $\epsilon >0$ is small, we conclude
$$\| f \| _{L^p_tW^{k,q}_x}\le  2 C(\omega )\delta +   2 C(\omega ) \epsilon O(\delta ^2 +\varepsilon
^2)+2C(\omega )C_1\varepsilon .
  $$
Then by a continuity argument, we conclude that (3.8) holds in $(0,T
_1)$ with $T_1=T$, i.e. the claim of Lemma 3.5.

\bigskip
Having obtained Lemma 3.5, we rewrite (3.9) in more precise form:

$$\aligned &if_t=\left \{ H_ {\omega _0}   + (\dot \gamma +\omega
-\omega _0) (P_+(\omega _0)-P_-(\omega _0))\right \}  f
\\& +
\sum _{  |m+n|=2} z^m \bar z^nP_c( \omega _{0})R_{m,n}(\omega )
  + \sum _{  |m + n|=1} z^m \bar z^n P_c( \omega _{0})A_{m,n}(\omega ) f+
\\& +
 (\dot \gamma +\omega -\omega _0) \left (P_c(\omega _0)\sigma _3-  (P_+(\omega _0)-P_-(\omega _0))  \right ) f  +
   \left (  V  (\omega ) -
 V  (\omega _0) \right )   f
\\& +(\dot \gamma +\omega -\omega _0) \left ( P_c(\omega  )-  P_c(\omega _0)\right )\sigma
_3f + O(f^2)+O_{loc}(z^3)\\& + ( P_c( \omega )-P_c( \omega _0))
\left (\sum _{ |m+n|=2} z^m \bar z^nR_{m,n}(\omega )
  + \sum _{  |m + n|=1} z^m \bar z^n  A_{m,n}(\omega ) f\right ) .
\endaligned \tag 3.11
$$
 We then set

$$ f_{2}=f+ \sum _{  |m+n|=2} R_{H_{\omega _0}}^{+}((m-n)\cdot \lambda (\omega _0)
) P_c( H _{\omega _0})R_{m,n}    (\omega  )z^m  \bar z^n.\tag 3.12
$$
We will need below:

\proclaim{Lemma 3.6} Assume the hypotheses of Lemma 3.5.  Then for
$s>1$ sufficiently large we can decompose $f _{2}= h_1+h_2+h_3+h_4$
with: {\item {(1)}} for a fixed $c_0(\omega _0), $ $\| h _1\|
_{L^2_t (\Bbb R, L^{2,-s}_x)} \le c_0(\omega _0) \| f(0)\| _{H^1}\le
c_0(\omega _0) \delta  ;$ {\item {(2)}} for a fixed $c_1(\omega _0),
$ $\| h _2\| _{L^2_t( \Bbb R, L^{2,-s}_x)} \le c_1(\omega _0)
|z(0)|^2\le c_1(\omega _0)\delta ^2;$ {\item {(3)}} $\| h _3\|
_{L^2_t( (0,T), L^{2,-s}_x)} \le O(\epsilon (\varepsilon + \delta )
) ;$ {\item {(4)}} for all admissible pairs $(r,p)$  we have $ \|
h_4 \| _{L^{r}_t((0,T), L^{p}_x)} =O(\epsilon \, \varepsilon  ).$

All the constants, included those in the big $O$'s, do not depend on
$T$.
\endproclaim
{\it Proof.} The proof is basically that in \S 4 \cite{CM}. We have
schematically

$$\aligned &i\partial _tP_c(H _{\omega _0})f_{2}=\left ( H _{\omega _0}   + (\dot \gamma +\omega -\omega _0)
(P_+(\omega _0) -P_-(\omega _0) )\right ) P_c(H _{\omega _0}) f_{2}
+
\\& + \sum _{|m|=2 } O(|z|^{3}) \left(
 R_{H _{\omega _0}}^+(m\cdot\lambda (\omega _{0} )  ) R_{m,0}
(\omega _{0} )+ R_{H _{\omega _0}}^+(-m\cdot\lambda (\omega _{0} ) )
R_{0,m}  (\omega _{0} )\right )
\\&  +
  P_c(H _{\omega _0})
  \left (o(1) O_{ loc }(|z|^{2}) +o(1) O_{ loc }(f)+O(f^2)\right )
   .  \endaligned  $$
For $h_1(0)=f(0)$ let
$$ i\partial _t (h_1+h_2)=\left ( H_{\omega _0}   + (\dot \gamma +\omega -\omega _0) (P_+-P_-)\right )
(h_1+h_2) , \quad  h_1(0)+h_2 (0)=f _{2}(0)  .$$ Then (1) follows by
Stricharz estimates applied to $P_\pm (\omega _0)h_1(t)=U_\pm (t,0)
f(0)$, with $U_\pm (t,s)$ defined in Lemma 3.5. To get (2) recall
from \cite{Cu3} that for a constant $C=C(\Lambda , \omega _{0})$
upper semicontinuous in $\omega _{0}$ and in $\Lambda
>\omega $ we have

$$
\|   U_\pm (t,t')R _{H_{\omega  }}^{+}(\Lambda     ) P_c g\|
_{L_x^{2,-s}} < C\langle t-t' \rangle ^{-\frac 32} \|  g \|
_{L_x^{2,s}} \, ,\, s> s _0 . \tag 5
$$
By  $f _{2}(0) =\sum _{ |m+n|=2} R_{H_{\omega _0}}^{+}((m-n)\lambda
(\omega _0 )) R_{m,n} (\omega _0)z^m(0)  \bar z^n (0)  $ and (5) we
get (2). Next we define   $h_3(0)=0$ and
$$\aligned &i\partial _t P_c(H _{\omega _0}) h_3=\left ( H_{\omega _0}   + (\dot \gamma +\omega -\omega _0)
(P_+(\omega _0) -P_-(\omega _0) )\right ) P_c(H _{\omega _0}) h_3 +
\\&
+ P_c(H _{\omega _0})\left (O(\epsilon ) O_{ loc }(|z|^{2})
+O(\epsilon ) O_{ loc }(f)+O(f^2)\right )   . \endaligned  $$  Then
(3) follows in a standard way from Strichartz inequalities, see
\cite{CM}. Finally we set $h_4(0)=0$ and

$$\aligned &i\partial _t h_{4}=\left ( H _{\omega _0}   + (\dot \gamma +\omega -\omega _0) (P_+-P_-)\right ) h_{4} + \\& +
\sum _{|m|=2 } O(|z|^{3}) \left(
 R_{H _{\omega _0}}^+(m\cdot\lambda (\omega _{0} )  ) R_{m,0}
(\omega _{0} )+ R_{H _{\omega _0}}^+(-m\cdot\lambda (\omega _{0} ) )
R_{0,m}  (\omega _{0} )\right ) .\endaligned $$ Then we have
$h_{4}=h_{41 }+h_{42}$ with $h_{4j}=\sum _\pm h_{4j\pm }$ with

$$h_{4 1\pm } (t)=\sum _{|m|=2 }\int _0^tU_\pm (t,t')
 O(|z(s)|^{3}) R_{H _{\omega _0}}^+(m\cdot\lambda (\omega _{0} )  ) R_{m,0}
(\omega _{0} ) dt' $$ and  $h_{42\pm }$ defined similarly. By (5) we
get $ \|h_{4j\pm } (t)  \|  _{L_x ^{2,-s}}\le C \epsilon \int _0^t
\langle t-t'\rangle ^{-\frac 32} |z(t')|^{2}  dt'$  and so  $ \| h
_4\| _{L^2_{t }L_x ^{2,-s}} \le \epsilon \| z \| _{L^{4} _{t }}
^{2}= O(\epsilon \, \varepsilon   ). $ This concludes the proof of
Lemma 3.6.

\bigskip

By substitution of (3.12) in the discrete part in  (3.7) we get
$$\aligned &i\dot z _j\xi _j-\lambda _j(\omega ) z_j\xi _j = P_{\ker (H  _\omega -\lambda _j)}
 ( \sum _{  |m+n|=2} z^m \bar
z^nR_{m,n}(\omega ) +\sum _{  |m+n|=3} z^m \bar
z^nR_{m,n}^{(1)}(\omega )\\& -\sum _{  |m' + n'|=1}\sum _{  |m +
n|=2} z^{m +m'}\bar z^{n+n'} A_{m',n'}(\omega ) R_{H_{\omega _0
}}^{+}((m-n)\cdot \lambda (\omega _0) ) P_c( H _{\omega })R_{m,n}
(\omega )\\&
 +\sum _{  |m + n|=1} z^m \bar z^n  A_{m,n}(\omega ) f_2+O(f^2)+O_{loc}(z^4).
\endaligned
$$
Here recall $P_{\ker (H  _\omega -\lambda _j)}= s _j \xi _j \langle
\quad , \sigma _3\xi _j \rangle $ with $s _j=\langle \xi _j , \sigma
_3\xi _j \rangle $ the signature of $\lambda _j$, that is either 1
or $-1$ and with $ s _1=-1$. By standard normal forms arguments
there exists a change of variables $ {\zeta }_j=z_j +\sum
_{|m+n|=2}^3p_{j,m,n}(\omega )z^m\overline{z}^n$ with $\Im
p_{j,m,n}=0$    for $|m+n|=2$ such that below we have $\Im a_{j,\ell
} (\omega )=0$ and

$$\aligned &i {\dot {\zeta}} _j\xi _j-\lambda _j(\omega ) \zeta_j\xi _j = \xi _j
\sum _{\ell}a_{j,\ell } (\omega )|\zeta _\ell |^2\zeta
_j+O_{loc}(\zeta f_2) +O(f^2)+O_{loc}(\zeta ^4)\\& - \sum _\ell
|\zeta _\ell |^2\zeta _jP_{\ker (H  _\omega -\lambda _j)}
A_{0,\delta _{\ell }}(\omega )R_{H_{\omega _0}}^{+}( \lambda _{\ell
}(\omega _0 ) +\lambda _j (\omega _0) ) P_c( H _{\omega })R_{\delta
_{\ell }+\delta _j,0} (\omega  )
\endaligned \tag  3.13
$$
where   $O_{loc}(\zeta ^n f_2)=\sum
 _\ell O_{loc}( \zeta _\ell ^{n}f_2)$.
Applying  $\langle \quad , \sigma _3\xi _j\rangle $ to (3.13)   and
recalling $P_{\ker (H  _\omega -\lambda _j)}=s _j\xi _j\langle \quad
, \sigma _3\xi _j\rangle $ we get
$$\aligned &( i {\dot {\zeta}} _j -\lambda _j(\omega ) \zeta_j) s _j = s _j
\sum _{\ell}a_{j,\ell } (\omega )|\zeta _\ell |^2\zeta _j+\langle
O_{loc}(\zeta f_2) +O(f^2)+O_{loc}(\zeta ^4) ,    \sigma _3\xi
_j\rangle  -\\&   -\sum _\ell |\zeta _\ell |^2\zeta _j s _j\langle
\xi _j , \sigma _3\xi _j\rangle
  \langle
A_{0,\delta _{\ell }}(\omega )R_{H_{\omega _0}}^{+}( \lambda _{\ell
}(\omega _0 ) +\lambda _j (\omega _0) ) P_c( H _{\omega })R_{\delta
_{\ell }+\delta _j,0} (\omega  ) , \sigma _3\xi _j\rangle .
\endaligned
$$
So
$$\aligned &  i {\dot {\zeta}} _j -\lambda _j(\omega ) \zeta_j  =
\sum _{\ell}a_{j,\ell } (\omega )|\zeta _\ell |^2\zeta _j+\langle
O_{loc}(\zeta f_2) +O(f^2)+O_{loc}(\zeta ^4) ,    \sigma _3\xi
_j\rangle  -\\& -s _j   \sum _\ell |\zeta _\ell |^2\zeta _j \langle
A_{0,\delta _{\ell }}(\omega )R_{H_{\omega _0}}^{+}( \lambda _{\ell
}(\omega _0 ) +\lambda _j (\omega _0) ) P_c( H _{\omega })R_{\delta
_{\ell }+\delta _j,0} (\omega  ) , \sigma _3\xi _j\rangle .
\endaligned \tag 3.14
$$
Recall that $\Im a_{j,\ell } (\omega )=0$. Then multiplying (3.14)
by $\overline{\zeta}_j$ we get
$$\aligned &   \frac{1}{2} \frac{d}{dt}{| \zeta _j|^2}   =
-s _j   \sum _\ell  \left ( \Gamma _{j,\ell } (\omega _0)
+o(1)\right ) |\zeta _\ell |^2|\zeta _j|^2+\\& +\Im [
\overline{\zeta _j}\langle O_{loc}(\zeta f_2) +O(f^2)+O_{loc}(\zeta
^4) ,    \sigma _3\xi _j\rangle ]
\endaligned \tag 3.15
$$
where we use $\omega -\omega _0=O(\epsilon )$,
$$\aligned & \Gamma _{j,\ell } (\omega  )=
 \Im \langle A_{0,\delta _{\ell
}}(\omega  )R_{H_{\omega   }}^{+}( \lambda _{\ell }(\omega   )
+\lambda _j (\omega  ) ) P_c( H _{\omega  })R_{\delta _{\ell
}+\delta _j,0} (\omega   ) , \sigma _3\xi _j(\omega   )\rangle
\endaligned $$
and the continuous dependence  in $\omega $   of $A_{0,\delta _{\ell
}}(\omega )$, $R_{\delta _{\ell }+\delta _j,0} (\omega  ) $ and
$P_c( H _{\omega })$. By Lemma 3.4 we have $\Gamma _{j,\ell }
(\omega  )\ge 0.$ Now we assume the following hypothesis:

\proclaim{Hypothesis 3.7} We suppose that  $\Gamma _{j,j } (\omega
)>0$    for any $j$ .
\endproclaim
{\it Remark. }Since $2\lambda _j(\omega  )  >\omega $  for any $j$,
Hypothesis 3.7 looks like a generic condition   by Lemma 3.4.

\bigskip
By $s _1=-1$, integrating (3.15)  we get for a fixed $\Gamma >0$ and
for $s\gg 1$ fixed
$$\aligned &   {| \zeta _1(t)|^2}  \ge  {| \zeta
_1(0)|^2} +\Gamma \int _0^t | \zeta _1(\tau )|^4 d\tau -c(\omega )
\left ( \int _0^t | \zeta  (\tau )|^4 d\tau  \right ) ^{\frac{1}{2}}
\| f_2\| _{L^2 ((0,t),L_x^{2,-s})}  \\& +  o(1 ) \left (\| f\| _{L^2
((0,t),L_x^{2,-s})}^2+   \int _0^t | \zeta  (\tau )|^4 d\tau \right
).
\endaligned \tag 3.16$$
Similarly, for all the $j$ we have
$$\aligned &
\Gamma \int _0^t | \zeta _j(\tau )|^4 d \tau \le {| \zeta _j (t)|^2}
+ {| \zeta  _j(0)|^2} +c(\omega ) \left (  \int _0^t | \zeta (\tau
)|^4 d\tau  \right ) ^{\frac{1}{2}} \| f_2\| _{L^2
((0,t),L_x^{2,-s})}
\\& + o(1 ) \left (\| f\| _{L^2 ((0,t),L_x^{2,-s})}^2+ \int _0^t |
\zeta (\tau )|^4 d\tau \right ).
\endaligned \tag 3.17$$
For any fixed $C_1\gtrsim 1$ there is an $T>0$ such that we have $\|
\zeta _j\| _{L^4(0,T)}^2\le C_1\epsilon $. By Lemma 3.5 we conclude
$\| f\| _{L^2_t((0,T), L^{2,-s}_x)}<C(C_1)\epsilon  $ by $\epsilon
\ge \delta
>0$  for $\epsilon $ sufficiently small.
By Lemma 3.6 we have $\| f_2\| _{L^2_t((0,T), L^{2,-s}_x)}<c_0\delta
+O(\epsilon ^2  ) $ for a fixed constant $c_0\approx 1$. Then  by
(3.17) and by orbital stability we get $\| \zeta  \|
_{L^4(0,T)}^2\le C_0\epsilon $  for some fixed $C_0>0$. By a
continuity argument, the same argument can be repeated for any
$T>0$. Then $\| \zeta  \| _{L^4(0,\infty )}^2\le C_0\epsilon $. This
and the fact that $\dot \zeta \in L^\infty _t$ implies $\lim _{t\to
\infty}\zeta (t)=0.$ We claim that this is incompatible with (3.16).
For $t$ sufficiently large $| \zeta
 (t)|^2<\delta ^2/2 $. Then by (3.17) we get for a fixed $c$

$$\aligned &
(\Gamma +o(1)) \int _0^t | \zeta  (s)|^4 ds \le  c\delta ^2
+c(\omega ) \left ( \int _0^t | \zeta  (s)|^4 ds \right )
^{\frac{1}{2}} \| f_2\| _{L^2 ((0,t),L_x^{2,-s})} \\& + o(1 )
 \| f\| _{L^2 ((0,t),L_x^{2,-s})}^2  .
\endaligned  \tag 3.18$$
If $\Gamma \left ( \int _0^t | \zeta  (s)|^4 ds \right )
^{\frac{1}{2}} \ge 2 c(\omega )\| f_2\| _{L^2 ((0,t),L_x^{2,-s})} $,
then by (3.18)
$$ \aligned & \int _0^t | \zeta  (s)|^4 ds\approx \left(   \Gamma /2 +o(1)   \right )
      \int _0^t | \zeta  (s)|^4 ds   \le
c\delta ^2 + o(1 )
 \| f\| _{L^2 ((0,t),L_x^{2,-s})}^2  . \endaligned
$$
Then, by the argument in Lemma 3.5 we get $ \| f\| _{L^2
((0,t),L_x^{2,-s})}\le c_1\delta $ for fixed $c_1$and this implies
$\int _0^t | \zeta (s)|^4 ds \le  c_2\delta ^2$ for fixed $c_2$. If
instead for $ \varepsilon =\left ( \int _0^t | \zeta  (s)|^4 ds
\right ) ^{\frac{1}{2}}$ we have $\Gamma \varepsilon < 2 c(\omega
)\| f_2\| _{L^2 ((0,t),L_x^{2,-s})} $, then by Lemma 3.6 we get
$\varepsilon \le c_0\delta +O(\delta ^2+\epsilon (\varepsilon +
\delta ))$. This implies $\varepsilon \le c_2\delta $ for fixed
$c_2$. So in all cases we get  $\| \zeta \| _{L^4(0,\infty )}^2\le
c_2\delta .$ We pick now initial conditions  $f(0)=0$,  $|\zeta
_1(0)|=\delta$ and $\zeta _j(0)=0$ for $j>1$. Then, by (3.16) we get

  $$\aligned &   -\delta ^2/2\ge    \Gamma
  \int _0^t | \zeta _1(\tau )|^4 d\tau -c(\omega ) \left (
\int _0^t | \zeta  (\tau )|^4 d\tau  \right ) ^{\frac{1}{2}} \|
f_2\| _{L^2 ((0,t),L_x^{2,-s})}  \\& +  o(1 ) \left (\| f\| _{L^2
((0,t),L_x^{2,-s})}^2+   \int _0^t | \zeta  (\tau )|^4 d\tau \right
).
\endaligned   $$
By Lemma 3.5 we conclude $\| f\| _{L^2 ((0,\infty ),
L^{2,-s}_x)}<C(c_1)\delta $. Since  $f(0)=0$,  by Lemma 3.6 we have
$\| f_2\| _{L^2 ((0,\infty ), L^{2,-s}_x)}<C \epsilon \delta
 $. So we have
 $$ -\delta ^2/2> \Gamma \| \zeta _1\| _{L^4(0,\infty
)}^4+ O(  o(1) \delta ^2 )\ge O(  o(1) \delta ^2 ) .\tag 3.19 $$ By
$|O(  o(1) \delta ^2 )|\ll \delta $ it follows that (3.19) is
absurd. But (3.19) is a consequence of the fact that we are assuming
 the orbital stability of $e^{it\omega _0} \phi
_{\omega _0}$. This implies that $e^{it\omega _0} \phi _{\omega _0}$
is orbitally unstable. Theorem 3.2 is proved.

\head \S 4  The case when the internal modes are not necessarily
close to the continuous spectrum
\endhead

In this section we we consider the following two hypotheses:

{\item {(H8')}} $H_\omega$ has a certain number of
  simple positive eigenvalues with $0<N_j\lambda _j(\omega )<
\omega < (N_j+1)\lambda _j(\omega )$ with $N_j$ integers with
$N_j\ge 1$. $H_\omega $ does not have other eigenvalues and $\pm
\omega $ are not resonances. We set $N=\max _j N_j$.

{\item {(H9')}} For   multi indexes  $m=(m_1,m_2,...)$ and
$n=(n_1,...)$, setting $\lambda (\omega )=(\lambda _1(\omega ),...)$
and $(m-n)\cdot \lambda =\sum (m_{j}-n_{j})  \lambda _j  $, we have
the following two non resonance hypotheses: {\item {(i)}}
$(m-n)\cdot \lambda (\omega )=0$ implies $m=n$ if $|m|\le N+2$ and
$|n|\le N+2$;  {\item {(ii)}} $(m-n)\cdot \lambda (\omega )\neq
\omega $ for all $(m,n)$ with $|m|+|n|\le N+2$.

{\item {(H10')}} We assume the non degeneracy  Hypotheses 4.4 and
4.5.

{\item {(H11')}} $\beta (t)\in C^{N+2}(\Bbb R, \Bbb R)$.

 Under hypothesis (H8'), if $\xi$ is a generator of $\ker (H_\omega
-\lambda _j)$ for any $j$, then $\langle \xi , \sigma _3\xi \rangle
\neq 0$.

\proclaim{Theorem 4.1 } Under hypotheses (H1-7) and (H8'-11')  the
excited states $e^{i\omega t}\phi _{\omega }(x)$ are orbitally
unstable.

\endproclaim
The structure of the proof is similar to Theorem 3.2, only more
complicate because we perform the normal form argument in \cite{CM}.
Given two vectors we will write $\overrightarrow{a}\le
\overrightarrow{b}$ if $a_j\le b_j$ for all components. If this
happens  we write $\overrightarrow{a}< \overrightarrow{b}$ if we
have $a_j<b_j$ for at least one $j$.
 We will  set $(m-n)\cdot \lambda =\sum _j (m-n)_j  \lambda
_j$. We will   say that $ {m}\in Res$  if:

 {\item{ (i) } } $ {m}=(m_1,m_2,...)$ with
 $m_j\in \Bbb N \cup \{  0\}$ for all $j$;

{\item{ (ii) } }  $
  {m}\cdot \lambda (\omega )
>\omega $;

{\item{ (iii) } } Given an $ {n}$ satisfying (i) and with $ {n}<
{m}$ with $ {m}\in Res$, then  $
  {n}\cdot \lambda (\omega )
<\omega $.
\bigskip
We assume by absurd that the excited state   $e^{i\omega _0 t} \phi
_ {\omega _0} (x)$ is orbitally stable. We rewrite formula (3.4) and
we consider the expansion

$$\aligned & \sigma _3 \sigma _1 G'(R)=\sum _{   |m+n|=2}^{ 2N+1}
 R_{m,n}(\omega ) z^m  \bar z^n+
\sum _{ |m + n|=1} ^N z^m  \bar z^n A_{m,n}(\omega ) f+
O(f^2)+\cdots .
\endaligned   $$
  $A_{m,n}(\omega )$ and $R_{m,n}(\omega )$ are real with
$\sigma _1 {R}_{m,n}  =-  {R}_{n,m}  $ and $A_{m,n}  =-\sigma _1
A_{n,m}  \sigma _1 $. We  express (3.4) as

$$\aligned &if_t=\left ( H _{ \omega (t)}+P_c(H_\omega )\sigma _3 \dot \gamma \right )f  +
\sum _{  |m+n|=2}^{N+1} z^m \bar z^nP_c(H_\omega
)R_{m,n}^{(1)}(\omega )
\\&  + \sum _{  |m + n|=1}^{N} z^m \bar z^n P_c(H_\omega )A_{m,n}^{(1)}(\omega ) f+
O(f^2)+\sum _{m\in Res }O_{loc}(|z ^{m}|),
\endaligned \tag  4.1
$$
 and  for $\widetilde{R}_{m,n}^{(1)}= {R}_{m,n}^{(1)}$
$$\aligned i\dot z _j\xi _j-\lambda _j(\omega ) z_j\xi _j&= P_{\ker (H  _\omega -\lambda _j)}
 ( \sum _{  |m+n|=2}^{2N+1} z^m \bar
z^n\widetilde{R}_{m,n}^{(1)}(\omega ) \\&
 +\sum _{  |m + n|=1}^{N} z^m
 \bar z^n A_{m,n}^{(1)}(\omega ) f+O(f^2)+
 \sum _{m\in Res }O_{loc}(|z ^{m}|)).
\endaligned \tag  4.1
$$
We consider $k=1,2,... N$ and set $f=f_k$   and $ z   _ {(k
),j}=z_j$ for $k=1$. The other $f_k$  and $  z   _ {(k),j}$ are
defined  below by induction. In (4.1) for $\ell =1$ we have
$$\aligned & \sigma _1f_\ell =\overline{f}_\ell , \text{ $A_{m,n}^{(\ell )} $,
  $R_{m,n}^{(\ell )} $
  and $\widetilde{R}_{m,n}^{(\ell )} $  are
   real},\\&  \text{exponentially decreasing in $x$ and
$C^1$ in $(\omega , x) $; $\sigma _1R_{m,n}^{(\ell )} =-
R_{n,m}^{(\ell )} $}.
\endaligned \tag 4.2$$

 In the ODE's there will be error terms
of the form
$$E_{ODE}(k)=\sum _{M\in Res}\left \{ O( |z ^{M}_{(k)}|^{2  }
 )+O(  z ^{M}_{(k)}   f _{k}  )\right \} +O(f^2_{k})+
O(\beta (|f _{k}|^2f _{k})).$$ In the PDE's there will be error
terms   of the form
$$ E_{PDE}(k)= \sum _{M\in Res} O_{loc}( |z_{(k)}|  ^{M}|
 )|z _{(k)}| +O_{loc}(  z _{(k)} f _{k}  )+O(f^2_{k})+
O(\beta (|f _{k}|^2f _{k})).$$   Then we define   $f_1=f$ and,
summing only over $(m,n)$  with $|(m-n)\cdot \lambda |<\omega $,
$$f_{k }=f_{k-1}+\sum _{   |m + n|= k
 } R_{H_{\omega }}((m-n)\cdot\lambda )
P_c(H_\omega )R_{m,n}^{(k-1)}   (\omega )z^m_{(k-1)} \bar
z^n_{(k-1)} . \tag 4.3$$ By
$\sigma_1R_{m,n}^{(k-1)}=-R_{n,m}^{(k-1)}$, by the fact that
$R_{m,n}^{(k-1)}$ is real and by $\sigma_1H_\omega =-H_\omega
\sigma_1$ we get $\sigma _1f_k=\overline{f}_k$. Starting from
$\widetilde{R}_{m,n}^{(1)} = {R}_{m,n}^{(1)}$ and summing only over
$(m,n)$ with $\lambda _j(\omega )\neq (m-n)\cdot \lambda (\omega )$,
we set
$$z_{(k),j}\xi _j=z_{(k-1),j}\xi _j+\sum _{    |m + n|= k }
 P_{\ker (H  _\omega -\lambda _j)}\widetilde{R}_{m,n}^{(k-1)}
 (\omega )\frac{z^m_{(k-1)} \bar
z^n_{(k-1)}}{\lambda _j-(m-n)\cdot \lambda }  .\tag 4.3
$$
We get the   equations

$$\aligned &i\partial _t f_{k}=\left ( H _{\omega }  +\sigma _3 \dot \gamma \right  )f_{k}
+ E_{PDE}(k)+\\&  \sum _{k+1\le |m +n|\le N +1 }
R_{m,n}^{(k)}(\omega ) z^m _{(k)}\bar z^n _{(k)}\text{ (sum over
pairs with }|(m-n)\cdot \lambda |<\omega )\\& +\sum _{2\le |m +n|\le
N+1} R_{m,n}^{(k)}(\omega ) z^m _{(k)}\bar z^n_{(k)} \text{ (sum
over pairs with }|(m-n)\cdot \lambda |>\omega ) ;
\endaligned \tag 4.4_k
$$
$$\aligned &
 i   {\dot z _{(k),j}}  \xi _j-\lambda _j(\omega ) z _{(k),j}\xi _j = \xi _j
 \sum _{   |m|=1} ^{N} a_{j,m  }^{(k)}(\omega )|z _{(k) } ^{m}| ^{2 }
 z _{(k),j}+   E_{ODE}(k) + \\& \sum _{  |m +n|= k +1 } ^{  2N +1 }
P_{\ker (H _\omega -\lambda _j)}\widetilde{R}_{m,n}^{(k)}(\omega )
z^m _{(k)}\bar z^n _{(k)}  \text{ ( with $ (m-n)\cdot \lambda
 \neq \lambda _j $)}\\& + \sum _{ |m +n|=1 }^{N} {z } ^{m }_{(k) }
  \overline{z  } ^{n }_{(k) }  P_{\ker (H  _\omega -\lambda _j)}
    {A}_{m,n }^{(k)}(\omega ) f _{k}.
\endaligned   \tag 4.4_k$$
The coefficients   in ($4.4_k$) are real because their entries are
products of entries of the coefficients   in (4.3),   which are
real, with Taylor coefficients of the rhs in ($4.4_{k-1}$)   at $z
_{(k-1)}=\bar z _{(k-1)}=0$ and $f_{k-1}=0$, which are also real. So
 in particular $\Im   [  a_{j,m
}^{(k)}(\omega )  ] =0 $. Since by (4.3) and by induction we have
$\sigma _1f_k=\overline{f}_k,$ taking complex conjugate in the $f_
k$ equation  in ($4.4_k$) we get $ \sigma _1R_{m,n}^{(k)} =-
R_{n,m}^{(k)} $.   At the step $k=N$, we can define

$$\aligned & \zeta _{ j}=z_{(N),j} + p_j
(z_{(N) },\overline{z}_{(N) })+ \sum _{1\le |m +n|\le N }z^m_{(N)}
\overline{z}^n_{(N)}\langle f_{N}, \alpha _{jmn} \rangle , \text{
with:}\endaligned$$    $\alpha _{jmn} $ vectors with entries which
are real valued exponentially decreasing functions; $p_j$
polynomials in $(z_{(N) },\overline{z}_{(N) })$ with real
coefficients and whose monomials have degree not smaller than $N+1$.
The above transformation can be chosen so that:

$$\aligned &
 i   {\dot \zeta _{ j}}  \xi _j-\lambda _j(\omega ) \zeta _{ j}\xi _j = \xi _j
 \sum _{ 1\le |m|\le N} b_{j,m  } (\omega )|\zeta ^{m}| ^{2 }
 \zeta +   E_{ODE}   +\\& + \sum _{ n+\delta _j\in Res}
  \overline{\zeta  } ^{n }
  P_{\ker (H  _\omega -\lambda _j)}
   { {A}}_{0,n }^{(N)} (\omega ) f _{N}
\endaligned \tag 4.5$$
with $   b_{j,m  } (\omega )  $ real and $E_{ODE}$ an error term
$$E_{ODE} =\sum _{M\in Res}\left \{ O( |\zeta ^{M} |^{2  }
 )+O(  \zeta ^{M}    f _{N}  )\right \} +O(f^2_{N})+
O(\beta (|f _{N}|^2f _{N})) .$$ We write

$$\aligned &i\partial _t P_c(\omega _0)f_{N}=\left \{ H_
{\omega _0}   + (\dot \gamma +\omega -\omega _0) (P_+(\omega
_0)-P_-(\omega _0))\right \} P_c(\omega _0)f_{N} +\\& +P_c(\omega
_0)\widetilde{E}_{PDE}(N) +  \sum _{2\le |m +n|\le N+1} P_c(\omega
_0)R_{m,n}^{(N)}(\omega  _0) \zeta ^m \bar \zeta ^n
\endaligned  $$ where $|(m-n)\cdot \lambda |>\omega$ for $|m +n|\le N$
and with

$$\aligned & \widetilde{E}_{PDE}(N)= E_{PDE}(N)  + \sum
_{2\le |m +n|\le N+1} P_c(\omega _0)\left ( R_{m,n}^{(N)}(\omega
 )-R_{m,n}^{(N)}(\omega _0) \right ) \zeta ^m \bar \zeta ^n
+\\&  +
 (\dot \gamma +\omega -\omega _0) \left (P_c(\omega _0)\sigma _3-  (P_+(\omega _0)-P_-(\omega _0))  \right ) f_N +
   \left (  V  (\omega ) -
 V  (\omega _0) \right )   f_N
\\& +(\dot \gamma +\omega -\omega _0) \left ( P_c(\omega  )-  P_c(\omega _0)\right )\sigma _3f_N.
  \endaligned
$$
Then proof of the following is almost the same of Lemma 3.5:

\proclaim{Lemma 4.2} Let $0<\delta < \epsilon $ be as in the
Definition
  3.1 of linear stability. For any $C_1>0$
there are a $\varepsilon (C_1)>0$ and a $C(C_1)$ such
  that if, for $0<\varepsilon <\varepsilon (C_1)$, we have $\sum _{m\in Res}\|
 { \zeta }^m \| _{L^{ 2 }_t(0,T)}  \le  C_1\epsilon $  for all $j$, then
for all admissible pairs $(p,q)$ we have for a fixed $c_0$
$$\| f_N\| _{L^p_t((0,T), W^{1,q}_x)}<c_0\delta +C(C_1)\varepsilon . $$
\endproclaim

  For $ |(m-n)\cdot \lambda (\omega  _{0})|>\omega _0$ if
  $|m+n|\le N$ in the sum (4.6) below, we set
$$\aligned &f_{N}=
  -\sum _{2\le |m+n|\le N+1}
  R_{H_ {\omega _0} }^+((m-n)\cdot \lambda (\omega  _{0})  )
 P_c(\omega
 )R_{m,n}^{(N)}   (\omega _{0} )\zeta ^m
  \bar \zeta ^n     + f_{N+1}. \endaligned \tag
4.6$$ Then  we get
$$\aligned &i\partial _t P_c(\omega _0) f_{N+1}=\left ( H _{\omega _0}   + (\dot \gamma +\omega -\omega _0)
(P_+(\omega _0) -P_-(\omega _0) )\right ) P_c(\omega _0) f_{N+1} +
\\&   \sum   O(|\zeta |^{|m +n|+1}
)R_{H _{\omega _0}}^+((m-n)\cdot\lambda (\omega _{0} )  )
R_{m,n}^{(N)} (\omega _{0} )
   +
  P_c(\omega _0)\widetilde{E}_{PDE}(N).  \endaligned \tag 4.7$$
We have $O(|\zeta |^{|m +n|+1} )=O(|\zeta ^M\zeta|)$ for $M\in Res$
for each factor in (4.7). By a proof similar to Lemma 3.6, see also
\cite{CM}, we have: \proclaim{Lemma 4.3} Assume the hypotheses of
Lemma 4.2. Then for $s>1$ sufficiently large we can decompose $f
_{N+1}= h_1+h_2+h_3+h_4$ with: {\item {(1)}} for a fixed $c_0(\omega
_0), $ $\| h _1\| _{L^2_t(\Bbb R, L^{2,-s}_x)} \le c_0(\omega _0) \|
f(0)\| _{H^1}\le c_0(\omega _0) \delta  ;$ {\item {(2)}} for a fixed
$c_1(\omega _0), $ $\| h _2\| _{L^2_t(\Bbb R, L^{2,-s}_x)} \le
c_1(\omega _0) |z(0)|^2\le c_1(\omega _0)\delta ^2;$ {\item {(3)}}
$\| h _3\| _{L^2_t((0,T), L^{2,-s}_x)} \le O(\epsilon (\varepsilon +
\delta ) ) ;$ {\item {(4)}} for all admissible pairs  $ \| h_4 \|
_{L^{r}_t ((0,T),L^{p}_x)} =O(\epsilon \, \varepsilon  ).$
\endproclaim

Substituting $ f_{N} $  in (4.5)    with the right hand side of
(4.6) we get

$$\aligned &
 i   {\dot \zeta _{ j}}  \xi _j-\lambda _j(\omega ) \zeta _{ j}\xi _j = \xi _j
 \sum _{ 1\le |m|\le N} b_{j,m  } (\omega )|\zeta ^{m}| ^{2 }
 \zeta
+\sum _{ 2\le |m+n|\le N+1}  \sum _{ n+\delta _j\in Res}
 \zeta ^m \overline{\zeta  } ^{\widetilde{n} +n}
 \times \\&  P_{\ker (H  _\omega -\lambda _j)}
   { {A}}_{0,n }^{(N)} (\omega ) R_{H_ {\omega _0} }^+((m-n)\cdot \lambda (\omega  _{0})  )
 P_c(\omega
 )R_{m,n}^{(N)}   (\omega _{0} ) \\&
+ \sum _{ n+\delta _j\in Res}
  \overline{\zeta  } ^{n }
  P_{\ker (H  _\omega -\lambda _j)}
   { {A}}_{0,n }^{(N)} (\omega ) f _{N+1}
+   E_{ODE}
\endaligned  $$
where $|(m-n)\cdot \lambda |>\omega$ for $|m +n|\le N$ in the above
formula. We
    considerate
new change of variables $\widetilde{\zeta }_j= \zeta _j+p _j(
{\zeta} ,\overline{\zeta})$ such that
$$\aligned &
 (i\dot { \widetilde{\zeta} _j}-\lambda _j(\omega )\widetilde{\zeta }_j )\xi _j=
\xi _j \sum _{ 1\le |m|\le N}  \widetilde{{a}}_{j,m  }(\omega
)|\widetilde{\zeta} ^{m}| ^{2 }
 \widetilde{\zeta} _j+   E_{ODE}(N)   - \\&     \sum _{m+\delta _j\in Res}
  |\widetilde{\zeta } ^m|^2\widetilde{\zeta }_j  P_{\ker (H  _\omega -\lambda _j)}
    {A}_{0,m}^{(N)}(\omega )
  R_{H _{\omega _0} }^+(m \cdot \lambda  (\omega _{0} )+ \lambda _j (\omega _{0} )  )
  R_{m+\delta _j,0}^{(N)}(\omega _{0} )
\\& - \sum _{m+\delta _j\in Res}\overline{\widetilde{\zeta}} ^m P_{\ker (H  _\omega -\lambda _j)}
  {A}_{0,m}^{(N)}(\omega ) f
_{N+1}
\endaligned \tag 4.8$$
with $\widetilde{a}_{j,m  }$, $ {A}_{0,m}^{(N)}$ and $R_{m+\delta
_j,0}^{(N)}$ real and with all the $m$ such that $m+\delta _j\in
Res$. It is possible to choose  $ p_j(\omega , z,\overline{z})$  as
polynomials with monomials $z ^m \overline{z} ^{n +\widetilde{n} }$
which, by $(m+ {n})\cdot \lambda >\omega$, are $O(z^{M})$ for $M\in
Res$. This implies $\sum _{M\in Res }\| \zeta ^M(t)\|_{ L_t ^{2
}}\approx \sum _{M\in Res }\| \widetilde{\zeta} ^M(t)\|_{ L_t ^{2
}}$.

Applying $\langle \quad , \sigma _3\xi _j\rangle $ to (4.8) and
recalling $P_{\ker (H  _\omega -\lambda _j)}=s _j\xi _j\langle \quad
, \sigma _3\xi _j\rangle $ we get
$$\aligned & i {\dot {\widetilde{\zeta}}} _j -\lambda _j(\omega ) \widetilde{\zeta}
_j  =   \sum _{ 1\le |m|\le N} \widetilde{{a}}_{j,m }(\omega
)|\widetilde{\zeta} ^{m}| ^{2 }
 \widetilde{\zeta} _j+\langle  E_{ODE}(N) \sigma _3\xi _j\rangle
 -s _j\times \\&     \sum _{m+\delta _j\in Res}
  |\widetilde{\zeta } ^m|^2\widetilde{\zeta }_j
  \langle    {A}_{0,m}^{(N)}(\omega )
  R_{H _{\omega _0} }^+(m \cdot \lambda  (\omega _{0} )+ \lambda _j (\omega _{0} )  )
  R_{m+\delta _j,0}^{(N)}(\omega _{0} ) , \sigma _3\xi _j\rangle \\&
  -   s_j  \sum _{m+\delta _j\in Res}
  |\widetilde{\zeta } ^m|^2\widetilde{\zeta }_j
  \langle    {A}_{0,m}^{(N)}(\omega )
  f_{N+1}, \sigma _3\xi _j\rangle .
\endaligned  \tag 4.9
$$
We can denote by $\Gamma _{m+\delta _j,j}(\omega ,\omega _0)$ the
quantity $\Gamma _{m+\delta _j,j}(\omega ,\omega _0 )=$
$$\aligned & \Im \left (\langle   {A}_{0,m}^{(N)}(\omega )
  R_{H _{\omega _0}}^+(m \cdot \lambda  (\omega  _{0})+ \lambda _j (\omega _{0} )  ) R_{m+\delta _j,0}^{(N)}(\omega _{0} )
   ,\sigma _3 \xi _j (\omega )\rangle
\right )\\& =\pi
  \langle   {A}_{0,m}^{(N)}(\omega )
  \delta (H _{\omega _0} -  m \cdot \lambda  (\omega _{0} )- \lambda _j (\omega  _{0}) )P_c(\omega
_0)R_{m+\delta _j,0}^{(N)}(\omega  _{0}), \sigma _3 \xi _j (\omega
)\rangle ,
 \endaligned  \tag 4.10$$ by
     $\frac{1}{x-i0}=PV\frac{1}{x}+i\pi \delta
_0(x)$,  \cite{Cu2} and Lemma 4.1 \cite{Cu3}  and which can be
proved as in Lemma 3.4.   We formulate the following hypothesis,
which is a conjecture:

\proclaim{Hypothesis  4.4 } We have the identity $\Gamma _{m+\delta
_j,j}(\omega ,\omega   )=$
$$\aligned & =({m_j+1})
 \langle
  \delta (H _{\omega  } -  m \cdot \lambda  (\omega   )- \lambda _j (\omega  ) )P_c(\omega
 )R_{m+\delta _j,0}^{(N)}(\omega   ), \sigma _3 R_{m+\delta _j,0}^{(N)}(\omega   )\rangle
.\endaligned $$
 \endproclaim
If Hypothesis  4.4  holds, then proceeding as in Lemma 3.4 we get
$\Gamma _{m+\delta _j,j}(\omega ,\omega   )\ge 0.$  We then assume
the following hypothesis:

\proclaim{Hypothesis  4.5 (non degeneracy hypothesis)} We have
$\Gamma _{m+\delta _j,j}(\omega ,\omega   )> 0$ for any $j$ and any
$m+\delta _j\in Res$.
\endproclaim
Recall that $\Im a_{j,m } (\omega )=0$. Then by (4.10) we get
$$\aligned &  \frac{d}{dt} \frac{|{\widetilde {\zeta}_j}|^2}{2} =
 -s _j\sum _{m+\delta _j\in Res}
  \Gamma _{m+\delta _j,j}(\omega ,\omega _0 )  |\widetilde {\zeta}  ^{m}\widetilde {\zeta} _j|^2  + \Im \left
( \langle E_{ODE}(N) , \sigma _3 \xi _j(\omega ) \rangle
\overline{\widetilde {\zeta}}_j\right )
\\&   +s _j \sum _{m+\delta _j\in Res}\Im \left ( \langle
 {A}_{0,m}^{(N)}(\omega )f _{N+1} , \sigma _3 \xi _j
(\omega )\rangle \overline{\widetilde {\zeta}} ^m
\overline{\widetilde {\zeta}} _j\right ) .
\endaligned
$$
By $s _1=-1$ we get for a fixed $\Gamma >0$
$$\aligned &   {| \widetilde{\zeta} _1(t)|^2}  \ge  {| \widetilde{\zeta}
_1(0)|^2} +\Gamma \int _0^t | \widetilde{\zeta }_1(\tau )|
^{2N_1+2}d\tau + o(1 ) (\| f_{N+1}\| _{L^2_tL_x^{2,-s}}^2+ \sum _{m
\in Res} \| \widetilde{\zeta} ^m \| _{L^2_t }^2)
\endaligned \tag 4.11$$
and for all the $m \in Res$
$$\aligned &
\Gamma \int _0^t | \widetilde{\zeta} ^m(\tau )|^2 d\tau  \le \sum _j
{| \zeta _j(t)|^2} \\& +  \sum _j{| \widetilde{\zeta} _j(0)|^2}  +
o(1 ) (\| f_{N+1}\| _{L^2_tL_x^{2,-s}}^2+ \sum _{m \in Res} \|
\widetilde{\zeta} ^m \| _{L^2_t }^2).
\endaligned  \tag 4.12$$
By the same argument in \S 3 we conclude  for a fixed $C_1$
$$ \sum _{m \in Res}\| \widetilde{\zeta} ^m \| ^{ 2}_{L^2(0,\infty
)}<C_1\epsilon .
$$ Then   $\lim _{t\to \infty}\widetilde{\zeta}
_j(t)=0.$ We claim that   this is incompatible with (4.11).
 For $t$
  large $| \widetilde{\zeta}
 (t)|^2<\delta ^2/2 $. Then by (4.12) we get for a fixed $c$

$$\aligned &
(\Gamma +o(1)) \sum _{m \in Res}\int _0^t | \widetilde{\zeta}
^m(\tau )|^2 d\tau  \le c\delta ^2+ \\& +c(\omega ) \left ( \sum _{m
\in Res}\int _0^t |\widetilde{\zeta} ^m (\tau )|^2 d\tau  \right )
^{\frac{1}{2}} \| f _{N+1}\| _{L^2 ((0,t),L_x^{2,-s})}
  + o(1 )
 \| f_N\| _{L^2 ((0,t),L_x^{2,-s})}^2  .
\endaligned  $$
As in \S 3 this yields    $\sum _{m \in Res}\| \widetilde{\zeta} ^m
\| _{L^2(0,\infty )}^2\le c_2\delta  $ for fixed $c_2$. We pick
initial conditions $f_N(0)=0$,  $|\widetilde{\zeta} _1(0)|=\delta$
and $\widetilde{\zeta} _j(0)=0$ for $j>1$. Then by (4.11) we get

  $$\aligned &   -\delta ^2/2\ge    \Gamma
  \int _0^t | \widetilde{\zeta }_1(\tau )| ^{2N_1+2}
   d\tau \\& -c(\omega ) \left (
\sum _{m \in Res}\int _0^t  |\widetilde{\zeta} ^m (\tau )|^2 d\tau
\right ) ^{\frac{1}{2}} \| f_{N+1}\| _{L^2 ((0,t),L_x^{2,-s})}  \\&
+  o(1 ) \left (\| f_N\| _{L^2 ((0,t),L_x^{2,-s})}^2+  \sum _{m \in
Res} \int _0^t |\widetilde{\zeta} ^m (\tau )|^2 d\tau \right ).
\endaligned   $$
By Lemma 4.2 we conclude $\| f_N\| _{L^2 ((0,\infty ),
L^{2,-s}_x)}<C(c_1)\delta $. Since  $f_N(0)=0$,  by Lemma 4.3 we
have $\| f_{N+1}\| _{L^2 ((0,\infty ), L^{2,-s}_x)}<C \epsilon
\delta
 $. So as in \S 3 we have
 $$ -\delta ^2/2> \Gamma \| \widetilde{\zeta} _1\| _{L^{2N_1+2}(0,\infty
)} ^{N_1+1} + O(  o(1) \delta ^2 )\ge O(  o(1) \delta ^2 )
$$ which is
absurd.  So Theorem 4.1 is proved.

\head \S 5 Completion of proof of inequality (2.4) \endhead

In this section we assume that $\omega$  is either exceptional of
second  or third type.  We consider   a factorization  $V=B^\ast A$
of  $H_{\omega }=\sigma _3(-\Delta +\omega )+V$,   with $A,B$ in
$C^2$, with real entries and exponentially decreasing.    As in \S 2
we  e consider a real matrix $U_1(x)\in C^\infty _0$ and we consider
the perturbation $ H _{\omega , \epsilon }= H_{\omega }+\epsilon
U_1$. We also consider a factorization  $U_1(x)=B^\ast _1 (x)A (x)$
with $B^\ast _1 (x)$    with real entries and in $C^2 _0$. We will
consider    $z\notin   [\omega ,+\infty )$ with $\Re z >0$ and in
some fixed small neighborhood of $\omega$. This $z$ can be expressed
as  $z=\omega -\zeta ^2$ with $\Re \zeta
>0$ close to 0. Then we
  the write $R_0(\zeta ):=R_{\sigma _3(-\Delta +\omega )}(z)$  with
integral  kernel
$$R_0(x,y, \zeta  ) :=\frac {\sigma _3}{4\pi |x-y|} \left [
\matrix e^{-\zeta |x-y|} & 0 \\ 0 & e^{-\sqrt{2\omega -\zeta ^2}
|x-y|}
\endmatrix \right ] .  \tag 5.1$$
Notice that (5.1) can be continued analytically in $\Re \zeta <0$,
but that this continuation does not represent the resolvent
$R_{\sigma _3(-\Delta +\omega )}(\omega -\zeta ^2)$. A Taylor
expansion at $\zeta =0$ of (5.1)  yields
$$\aligned & R_0(x,y,
 \zeta  )=\frac {\sigma _3}{4\pi |x-y|} \left [ \matrix 1 & 0 \\ 0 &
e^{-\sqrt{2\omega} |x-y|}
\endmatrix \right ] -
\left [ \matrix \frac 1{4\pi } & 0 \\ 0 & 0
\endmatrix \right ] \zeta +O(\zeta ^2 ).
\endaligned$$
We consider the corresponding expansion of $R_0( \zeta  ) \in \Cal B
( H^{-1}_{s},  H^1_{-s})$, $s>5/2$,
$$R_0( \zeta  ) =R_0(0 )-\zeta G_1  +O(\zeta ^2 ).$$
We set
   $ Q(z) = A   R_{H_\omega }(z)   B^\ast _1$ and $
Q_{1,\epsilon }(z) = A  R_{H _{\omega , \epsilon } }(z )   B^\ast
_1= (1+ \epsilon   Q(z) )^{-1}  Q(z) .$ Recall that $\Cal V$ was the
vector space formed by eigenvalues and resonant vectors of $H_\omega
$ at $\omega$. In \cite{CP}, following \cite{JK}, it is proved:

\proclaim{Lemma 5.1} There is an $(1+AR_0(0) B^\ast ) $ invariant
splitting
$$L^2= \ker (1+AR_0(0 ) B^\ast  )
\oplus \ker ^\perp (1+ B R_0(0 ) A^\ast  ) . \tag 5.2$$ Let $s>1/2$.
Then $\Cal V=\ker   (1+R_0(0) B^\ast A )$. The map $\psi \to   \Psi
=-A\psi   $ is an isomorphism
$$\ker   (1+R_0(0) B^\ast A )  \subset L^2_{-s}
\to \ker (1+ A R_0(0) B^\ast  )  \subset L^2 .\tag 5.3$$ The inverse
map is $\Psi \to R_0(0) B^\ast \Psi .$

\endproclaim
We will denote by $P\oplus Q$ the projections  associated to (5.2).

\bigskip

We will suppose now that (2.4) is not true. This implies that $$\inf
\{ \langle \sigma _3H_\omega u, u \rangle : \| u\|_2=1, \, u\in
H^1\cap L_c^2(H_{\omega })\} =-\lambda < 0.$$ In particular, there
exists $u\in L_c^2(H_{\omega })$ unitary with $\sigma _3H_\omega
u=-\lambda u$. Then, by standard theory $u\in C^2$ with
$|u(x)|\lesssim e^{-\sqrt{\omega } |x|}.$

We distinguish now between the  cases when $\omega$ is   of second
and of third type.

\head \S 5.1 $\omega$  exceptional of second type \endhead

We assume here that $\omega$ is an eigenvalue but not a resonance.
 Let $P_0$ be the natural spectral
projection in $L^2$ on $\Cal V:=\ker H_\omega$. By   Corollary 4.4
\cite{CP} for $s>5/2$  and for $\zeta $ near 0  we have   in $   B (
H^{-1}_{s}, H^1_{-s})$ the expansion
$$R_{H_\omega }(\omega-\zeta ^2)= \zeta ^{-2} P_0
- \zeta ^{-1} P_0 V G_3V P_0 + O(1) $$ with $ G_3(x,y)= \frac
{1}{24\pi }\text{diag}( 1 , 0) |x-y|^2,$ i.e. the diagonal 2x2
matrix with (1,0) on the diagonal. We can write
$$Q(\omega-\zeta ^2)=
 \zeta ^{-2}A  P_0B^\ast _1
- \zeta ^{-1} A P_0 V G_3V P_0B^\ast _1  +Q_c (\omega -\zeta ^2),
$$ where $Q_c (\omega -\zeta ^2) $ admits an analytic extension for
$\zeta$ around $0$. We write

$$\aligned & Q_{1,\epsilon }(z) =
\left [ 1+\epsilon  (1+  \epsilon   Q_c(z) )^{-1}\zeta ^{-2}
K(\epsilon ,\zeta )\right ] ^{-1}
  (1+  \epsilon  Q_c(z) )^{-1}
 Q(z)  \\& K(\epsilon ,\zeta ):=A  P_0B^\ast _1 - \zeta   A P_0 V G_3V P_0B^\ast
_1.\endaligned$$ By the fact that $\omega $ is of positive
signature, there is a basis $\psi _j$ of $\ker H_\omega $ such that
$ \langle \psi _j, \sigma _3 \psi _k \rangle = \delta _{j,k}$ We can
pick $U_1$  such that  we also have $\langle  \sigma _3U_1  \psi _j,
\psi _k \rangle = \delta _{j,k} d_j$,  with $d_j\neq d_k$ for $j\neq
k$ and $d_j<0$ for all $j$. Since $K(\epsilon ,\zeta )$ is of rank
$\dim \ker (H_\omega )$, we can consider the equation
$$ \aligned & \text{det} \left [ \zeta ^{2}+\epsilon  (1+  \epsilon   Q_c(z)
)^{-1} K(\epsilon ,\zeta )\right ]=  \text{det} \left [ \zeta
^{2}+\epsilon   K(\epsilon ,\zeta )(1+  \epsilon   Q_c(z)
)^{-1}\right ]= \\& \text{det} \left [ \zeta ^{2}+\epsilon    \left
[  A \psi _j\frac{\langle  B^\ast _1\cdot ,\sigma _3 \psi _k
\rangle }{\langle \psi _j,\sigma _3 \psi _j \rangle } - \zeta A \psi
_j \frac{\langle  V G_3V P_0B^\ast _1 \cdot ,\sigma _3 \psi _k
\rangle }{\langle \psi _j,\sigma _3 \psi _j \rangle }
 \right ] +O(\epsilon ^2)\right ] \\& = \text{det} \left [ \zeta
 ^{2}+\epsilon
  \delta _{j,k} d_j +  \epsilon \zeta  \frac{\langle  V G_3V P_0U_1\psi _j ,\sigma _3 \psi _k \rangle }{\langle \psi _j,\sigma _3
\psi _j \rangle } +O(\epsilon  ^{2})\right ]=0.
\endaligned \tag 5.4$$
We consider for  values $\epsilon >0$ the $2\dim \Cal V$  solutions
$  \pm \sqrt{\epsilon |d_j|}+ O(\epsilon ^{\frac{3}{4}}) $. The
solutions $\zeta _j (\epsilon )= \sqrt{\epsilon |d_j|}+ O(\epsilon
^{\frac{3}{4}}) $ with $\Re \zeta _j (\epsilon )>0$ yield a number
of $ \dim \Cal V$ of distinct eigenvalues of $H_{\omega ,\epsilon}$
given by $z_j( \epsilon )= \omega -\zeta _j ^2(\epsilon )$  and with
$\dim \ker (H_{\omega ,\epsilon}-z_j( \epsilon ) )=1$. The roots
$\zeta $ of (5.4)  with $\Re \zeta <0$ give singularities of the
analytic continuation of $Q_{1,\epsilon }(\omega -\zeta ^2)$ which
do not correspond to eigenvalues of $H_{\omega ,\epsilon}$. By the
symmetry of $\sigma (H_{\omega ,\epsilon})$ with respect to the
coordinate axes, all the $z_j( \epsilon )$ are on $\Bbb R$. We claim
now that it is possible to choose generators $\psi _j(\epsilon )\in
\ker (H_{\omega ,\epsilon}-z_j( \epsilon ) )$ such that
$$\langle u, \sigma _3\psi _j(\epsilon )\rangle =o(1) \text{ for $\epsilon \searrow
0 $}.\tag 5.5$$ If (5.5) is true, then (2.4) is true by an argument
similar to the one in \S 2. So now we focus on (5.5).  We consider a
nonzero solution
$$ (1+AR_0( \zeta  _1(\epsilon )) (B^\ast +
\epsilon  B^\ast _1) ) \Psi _j (\epsilon ) =0.\tag 5.6$$ Then $\psi
 _j(\epsilon )     =R_0( \zeta  _1(\epsilon ))(B^\ast + \epsilon
B^\ast _1) \Psi  _j(\epsilon )$ is a nonzero element in $\ker
(H_{\omega ,\epsilon}-z_j( \epsilon ) )$. We consider equation (5.5)
using the splitting (5.2). Notice that $G_1V\psi =0$ for any $\psi
\in \ker H_\omega$. As a consequence
$$ \aligned &Q(1+AR_0( \zeta  _1(\epsilon ))B^\ast )P=\\&   Q  A( R_0( \zeta  _1(\epsilon ))-R_0(0))B^\ast
P=-\zeta _1QAG_1B^\ast Q+O(\epsilon )= O(\epsilon ).
\endaligned$$
Similarly $P(1+AR_0( \zeta  _1(\epsilon ))B^\ast )Q=  O(\epsilon )$.
Set $B^\ast(\epsilon )=B^\ast + \epsilon  B^\ast _1.$ Then from

$$  \left [ \matrix P(1+AR_0( \zeta  _1(\epsilon )) B^\ast(\epsilon )  )P & O(\epsilon  )\\  O(\epsilon  )  &
Q( 1+AR_0( \zeta  _1(\epsilon )) B^\ast(\epsilon )  )  Q
\endmatrix \right ]  \left [ \matrix  P\Psi  (\epsilon ) \\ Q\Psi  (\epsilon ) \endmatrix \right ] =0$$
we  get
$$Q\Psi  (\epsilon )=\left [ Q( 1+AR_0( 0) B^\ast +o(1)  )  Q \right ] ^{-1}O(\epsilon  )P\Psi  (\epsilon
)$$ and so $\| Q\Psi  (\epsilon ) \| _2\le C \epsilon \| P\Psi
(\epsilon ) \| _2$. Normalizing $\| P\Psi (\epsilon ) \| _2=1$ we
see that
$$ \aligned & \langle u , \sigma _3\psi _j(\epsilon )  \rangle =
 \langle u , \sigma _3R_0( \zeta  _1(\epsilon )) B^\ast(\epsilon )P\Psi _j(\epsilon )  \rangle
 +
   \langle u , \sigma _3R_0( \zeta  _1(\epsilon ))B^\ast(\epsilon )Q\Psi _j(\epsilon )  \rangle   .\endaligned  $$
We have $ |\langle u , \sigma _3R_0( \zeta  _1(\epsilon
))B^\ast(\epsilon )Q\Psi _j(\epsilon )  \rangle | \le C \| u\|
_{L^{2,s}}   \| Q\Psi  (\epsilon ) \| _2=O(\epsilon )$, where we use
$|u(x)|\lesssim e^{-\sqrt{\omega } |x|}   $ and $\langle x \rangle
^{N} |B^\ast(\epsilon )(x)|< C_N$ for all $x\Bbb R^3$ and $\epsilon
$ small. Similarly
$$\langle u , \sigma _3R_0( \zeta  _1(\epsilon )) B^\ast(\epsilon )P\Psi _j(\epsilon )
\rangle = \langle u , \sigma _3R_0( 0) B^\ast P\Psi _j(\epsilon )
\rangle +o(1)=o(1)$$ by $\sigma _3R_0( 0) B^\ast P\Psi _j(\epsilon
)\in N_g(H_\omega ^\ast )$ and $u\in N_g^\perp (H_\omega ^\ast )$.

\head \S 5.2 $\omega$  exceptional of third type \endhead

We assume here that $\omega$ is an eigenvalue and a resonance. In
particular we pick an appropriately normalized $\psi $ resonant
vector, see Lemma 3.1 \cite{CP}.  Then for $P_0$ and $G_3$ as in \S
5.1 we have
$$R_{H_\omega }(\omega-\zeta ^2)= \zeta ^{-2} P_0
- \zeta ^{-1} P_0 V G_3V P_0 + \zeta ^{-1} \psi \langle  \quad ,
\sigma _3 \psi \rangle + O(1)
$$ with $ G_3(x,y)= \frac {1}{24\pi }\text{diag}( 1 , 0) |x-y|^2,$
i.e. the diagonal 2x2 matrix with (1,0) on the diagonal. We can
write
$$Q(\omega-\zeta ^2)=
 \zeta ^{-2}A  P_0B^\ast _1
- \zeta ^{-1} A P_0 V G_3V P_0B^\ast _1 +\zeta ^{-1}A\psi \langle
B^\ast \quad , \sigma_3\psi \rangle +Q_c (\omega -\zeta ^2),
$$  where $Q_c (\omega -\zeta ^2) $ admits an analytic extension for
$\zeta$ around $0$. We write

$$\aligned & Q_{1,\epsilon }(z) =
\left [ 1+\epsilon  (1+  \epsilon   Q_c(z) )^{-1}\zeta ^{-2}
K(\epsilon ,\zeta )\right ] ^{-1}
  (1+  \epsilon  Q_c(z) )^{-1}
 Q(z)  \\& K(\epsilon ,\zeta ):=A  P_0B^\ast _1 - \zeta   A P_0 V G_3V P_0B^\ast
_1+  \zeta  A\psi \langle B^\ast _1 \quad , \sigma_3\psi \rangle
.\endaligned$$ We consider   a basis $\psi _j$ of $\ker H_\omega $
such that $ \langle \psi _j, \sigma _3 \psi _k \rangle = \delta
_{j,k}$ and $\langle  \sigma _3U_1 \psi _j, \psi _k \rangle = \delta
_{j,k} d_j$,  with $d_j\neq d_k$ for $j\neq k$ and $d_j<0$ for all
$j$. We can also add that $\langle  \sigma _3U_1 \psi _j, \psi
\rangle = 0$ for all $j$ and $\langle  \sigma _3U_1 \psi  , \psi
\rangle = d>0 $ (it is easy to see that there is a $U_1$ satisfying
all the above hypotheses). Since $K(\epsilon ,\zeta )$ is of rank
$1+\dim \ker (H_\omega )$ we consider the equation
$$ \aligned &
\text{det} \left [ 1+\epsilon \zeta ^{-2}   K(\epsilon ,\zeta )(1+
\epsilon   Q_c(z) )^{-1} \right ]=0.\endaligned $$ This means we are
considering the determinant of a matrix of the form $O(\epsilon
^2)+$

$$\aligned &   \left [ \matrix \zeta ^{2}+\epsilon    \left (
 \frac{\langle U _1\psi _j ,\sigma _3 \psi _k  \rangle
}{\langle \psi _j,\sigma _3 \psi _j \rangle } - \zeta \frac{\langle
V G_3V P_0U _1\psi _j  ,\sigma _3 \psi _k \rangle }{\langle \psi
_j,\sigma _3 \psi _j \rangle }
 \right )  &    &    - \zeta \epsilon \frac{\langle
V G_3V P_0U _1\psi   ,\sigma _3 \psi _k \rangle }{\langle \psi
_j,\sigma _3 \psi _j \rangle }   \\ 0 & & \zeta +\epsilon \langle
U_1 \psi , \sigma _3 \psi \rangle
\endmatrix \right ] ,
\endaligned $$ that is

$$\aligned &
\text{det} \left ( \left [ \matrix \zeta
 ^{2}\delta _{j,k} +\epsilon
  \delta _{j,k} d_j    &  0 \\  0 & \zeta + \epsilon d\endmatrix  \right ]  + O(\epsilon \zeta +\epsilon ^2)\right ) =0.
\endaligned  $$
Since $d>0$, once again there are  only $   \dim \ker (H_\omega )$
roots  $\zeta _j (\epsilon )= \sqrt{\epsilon |d_j|}+ O(\epsilon
^{\frac{3}{4}}) $ with $\Re \zeta _j (\epsilon )>0$ which yield $
\dim   \ker (H_\omega )$   distinct eigenvalues of $H_{\omega
,\epsilon}$ given by $z_j( \epsilon )= \omega -\zeta _j ^2(\epsilon
)$. We have   $\dim \ker (H_{\omega ,\epsilon}-z_j( \epsilon ) )=1$
and $z_j( \epsilon )\in \Bbb R$ and  $z_j( \epsilon )<\omega .$ By
the same argument of \S 5.1  it is possible to choose generators
$\psi _j(\epsilon )\in \ker (H_{\omega ,\epsilon}-z_j( \epsilon ) )$
such that (5.5) holds. Then (2.4) is true by an argument similar to
the one in \S 2.

Thus the proof of (2.4) is completed


\Refs\widestnumber\key{1997shire}



\ref\key{AHS} \by S.Agmon,  I.Herbst, E.Skibsted \paper Perturbation
of embedded eigenvalues in the generalized $N $-body problem \jour
Comm. Math. Phys. \vol 122 \yr 1989 \pages  411--438
\endref




\ref\key{BP}  \by V.Buslaev, G.Perelman
 \paper
On the stability of solitary waves for nonlinear Schr\"odinger
equations \inbook Nonlinear evolution equations\eds N.N. Uraltseva
\pages 75--98 \bookinfo Transl. Ser. 2, 164 \publ Amer. Math. Soc.
\yr 1995 \publaddr Providence, RI
\endref





\ref\key{BS} \by V.S.Buslaev, C.Sulem \paper On the asymptotic
stability of solitary waves of Nonlinear Schr\"odinger equations
\jour Ann. Inst. H. Poincar\'e. An. Nonlin.  \vol 20 \yr 2003 \pages
419--475
\endref


\ref\key{CL} \by T.Cazenave, P.L.Lions \paper Orbital stability of
standing waves for  nonlinear Schr\"odinger equations \jour Comm.
Math. Phys.  \vol 85 \yr 1982 \pages 549--561
\endref




\ref\key{CHM} \by C.Cruz-Sampedro,  I.Herbst, R.Martinez-Avendano
\paper Perturbations of the  Wigner-Von Neuman  potential leaving
the embedded eigenvalue fixed \jour  Ann H. Poincare\vol 3 \yr 2002
\pages 331--346
\endref



\ref\key{CoP} \by A.Comech, D.Pelinovsky \paper Purely nonlinear
instability of standing waves with minimal energy\jour  Comm. Pure
Appl. Math. \vol    56 \yr 2003  \pages  1565--1607
\endref


\ref \key{Cu1} \by S.Cuccagna \paper On asymptotic stability in
energy space of  ground states of NLS in 1D \paperinfo
 http://arxiv.org/
\endref




\ref \key{Cu2} \bysame \paper Stabilization of solutions to
nonlinear Schr\"odinger equations \jour Comm. Pure App. Math. \vol
54 \yr 2001 \pages 1110--1145
\endref

\ref\key{Cu3} \bysame \paper On asymptotic stability of ground
states of NLS\jour Rev. Math. Phys. \vol 15 \yr 2003 \pages 877--903
\endref





\ref \key{CM} \by S.Cuccagna, T.Mizumachi\paper On asymptotic
stability in energy space of ground states for nonlinear
Schr\"odinger equations \paperinfo
 http://arxiv.org/
\endref


\ref\key{CP} \by  S.Cuccagna, D.Pelinovsky\paper Bifurcations from
the endpoints of the essential spectrum in the linearized nonlinear
Schr\"odinger problem\jour  J. Math. Phys. \vol {46} \yr 2005 \pages
053520
\endref

\ref\key{CPV}\by S.Cuccagna, D.Pelinovsky, V.Vougalter \paper
Spectra of positive and negative energies in the linearization of
the NLS problem\jour Comm.  Pure Appl. Math. \vol 58 \yr 2005 \pages
1--29
\endref

\ref \key{CT} \by S.Cuccagna, M.Tarulli\paper On asymptotic
stability in energy space of  ground states of NLS in 2D \paperinfo
http://arxiv.org/
\endref






 \ref\key{JK} \by A.Jensen, T.Kato \paper Spectral properties of
 Schr\"odinger operators and time decay of the wave functions  \jour
 Duke Math. J. \vol 46 \yr 1979 \pages 583--611
 \endref



\ref\key{J} \by C.K.R.T.Jones \paper  An instability mechanism for
radially symmetric standing waves of  nonlinear Schr\"odinger
equations\jour Jour. Diff. Eq. \vol 71 \yr 1988 \pages 34--62
\endref






\ref\key{Gz}\by Zhou Gang\paper Perturbation Expansion and N-th
Order Fermi Golden Rule of the Nonlinear Schr\"odinger Equations
\jour  J. Math. Phys. \vol {48} \yr 2007 \pages 053509 \endref





\ref\key{GS}\by Zhou Gang, I.M.Sigal\paper
 Relaxation of Solitons in Nonlinear Schr\"odinger Equations with Potential
\paperinfo http://arxiv.org/abs/math-ph/0603060 \endref




\ref\key{Gr1} \by M.Grillakis \paper Analysis of the linearization
around a critical point of an infinite dimensional Hamiltonian
system \jour Comm. Pure Appl. Math  \vol 43 \yr 1990 \pages 299--333
\endref

\ref\key{Gr2} \bysame \paper Linearized instability for nonlinear
Schr\"odinger and Klein Gordon equations \jour Comm. Pure Appl.
Math.  \vol 41 \yr 1988 \pages 747--774
\endref





\ref\key{GSS1} \by M.Grillakis, J.Shatah, W.Strauss \paper Stability
of solitary waves in the presence of symmetries, I \jour Jour.
Funct. An.  \vol 74 \yr 1987 \pages 160--197
\endref

\ref\key{GSS2} \bysame \paper Stability of solitary waves in the
presence of symmetries, II \jour Jour. Funct. An.  \vol 94 \yr 1990
\pages 308--348
\endref









\ref\key{M1} \by T.Mizumachi \paper A remark on linearly unstable
standing wave solutions to NLS\jour   Nonlinear Analysis \vol 64
\yr 2006 \pages 657--676
\endref

\ref\key{M2} \bysame \paper Vortex solitons for 2D focusing
Nonlinear Schr\"odinger Equations \jour   Diff. Integral Equations
\vol 18 \yr 2005 \pages 431--450
\endref


\ref\key{M3} \bysame \paper Instability of bound states  for 2D
Nonlinear Schr\"odinger Equations \jour   Discr. Cont. Dyn. Systems
\vol 13 \yr 2005 \pages 413--428
\endref


\ref\key{M4} \bysame \paper Instability of vortex solitons for 2D
focusing NLS
\endref









\ref\key{Si} \by I.M.Sigal \paper Nonlinear wave and Schr\"odinger
equations. I.
 Instability of periodic and quasi- periodic solutions
\jour Comm. Math. Phys. \vol 153 \yr 1993 \pages 297--320
\endref




\ref\key{SW1} \by A.Soffer, M.Weinstein \paper Selection of the
ground state for nonlinear Schr\"odinger equations
 \jour Rev. Math. Phys. \vol 16 \yr 2004 \pages
977--1071
\endref





\ref\key{SW2}\bysame
 \paper
 Resonances, radiation damping and instability
in Hamiltonian nonlinear wave equations \jour Invent. Math. \vol 136
\yr 1999 \pages 9--74
\endref







\ref\key{T} \by T.P.Tsai \paper  Asymptotic dynamics of nonlinear
Schr\"odinger equations with many bound states\jour   J. Diff. Eq.
\vol  192  \yr 2003 \pages  225--282
\endref





\ref\key{TY1} \by T.P.Tsai, H.T.Yau \paper Asymptotic dynamics of
nonlinear Schr\"odinger equations: resonance dominated and radiation
dominated solutions\jour  Comm. Pure Appl. Math. \vol  55  \yr 2002
\pages 153--216
\endref


\ref\key{TY2} \bysame \paper Relaxation of excited states in
nonlinear Schr\"odinger equations \jour   Int. Math. Res. Not. \vol
31  \yr 2002 \pages 1629--1673
\endref


\ref\key{TY3} \bysame \paper Classification of asymptotic profiles
for nonlinear Schr\"odinger equations with small initial data \jour
Adv. Theor. Math. Phys. \vol  6  \yr 2002 \pages  107--139
\endref



\ref \key {TY4} \bysame  \paper  Stable directions for excited
states of nonlinear Schr\"odinger equations \jour Comm. Partial
Diff. Eq. \vol 27  \yr 2002 \pages  2363--2402
\endref




\ref\key{We1} \by M.Weinstein \paper Lyapunov stability of ground
states of nonlinear dispersive equations \jour Comm. Pure Appl.
Math.  \vol 39 \yr 1986 \pages 51--68
\endref


\ref\key{We2} \bysame  \paper Modulation stability of ground states
of nonlinear Schr\"odinger equations \jour Siam J. Math. Anal. \vol
16 \yr 1985 \pages 472--491
\endref








\endRefs
\enddocument